\documentclass[11pt]{article}

\usepackage{amssymb,amsmath,bm}
\usepackage{amsthm}

\usepackage{textcomp}
\usepackage{enumerate}      
\usepackage{graphicx}        
\usepackage{caption}

\usepackage{mathrsfs}

\usepackage{url}

\renewcommand{\setminus}{{\smallsetminus}}

\newcommand{\cp}[1]{\vcenter{\hbox{#1}}}

\usepackage{amssymb,amsmath,bm}
\usepackage{amsthm}
\usepackage{hyperref}
\usepackage{mathrsfs}
\usepackage{textcomp}
\usepackage{enumerate}
\usepackage{graphicx}
\usepackage{tikz}
\usepackage{verbatim}

\newtheorem{theorem}{Theorem}[section]
\newtheorem{lemma}[theorem]{Lemma}

\newtheorem{definition}[theorem]{Definition}

\newtheorem{conjecture}[theorem]{Conjecture}

\theoremstyle{remark}
\newtheorem{remark}[theorem]{Remark}

\theoremstyle{remark}

\numberwithin{equation}{section}

\begin{document}
\title{\bf Computation of leading coefficients in asymptotics of\\ relative quantum invariants}

\author{Ka Ho Wong and Tian Yang}

\date{}

\maketitle

\begin{abstract} 

We propose Asymptotic Expansion Conjectures of the relative Reshetikhin-Turaev invariants, of the relative Turaev-Viro invariants and of the discrete Fourier transforms of the quantum $6j$-symbols Our conjectures refine the corresponding volume conjectures by identifying not only the exponential terms with volume, but also the leading coefficients in the asymptotics expansions in terms of adjoint twisted Reidemeister torsions, logarithmic holonomies, determinant of Gram matrices and edge lengths. For families of special cases for which the exponential terms are already known, we prove our conjecture by computing the leading coefficients in the asymptotic expansions. The significance of these expansions is that we do not specify the way that the sequence of the colorings converges to the limit. As a consequence, the terms in the expansion will have to depend on the index $r,$ but the dependence is in a way that  the terms  are purely geometric invariants  of the metrics on the underlying manifold and only the metrics vary with $r.$  Also, the appearance of the logarithmic holonomy and edge length terms in the relative setting is a new phenomenon that has never been observed in the expansion of either the colored Jones polynomials of links or the Reshetikhin-Turaev invariants of closed $3$-manifolds. 
\end{abstract}

\section{Introduction}

Volume Conjectures relate the exponential growth rate of various quantum invariants of hyperbolic $3$-manifolds to the hyperbolic volume of the manifolds. Examples include the relationships between: the values of the colored Jones polynomials of a hyperbolic  link at certain primitive roots of unity and the complete hyperbolic volume of the link complement\,\cite{Ka2, MM} and the values of the colored Jones polynomials  of a link near the roots of unity and the volume of the incomplete hyperbolic metrics on the link complement\,\cite{G}, the Reshetikhin-Turaev and the Turaev-Viro invariants of a hyperbolic $3$-manifold and the hyperbolic volume of the manifold\,\cite{CY}, the relative Reshetikhin-Turaev invariants of a closed oriented $3$-manifold with a hyperbolic link inside it and the volume of the hyperbolic cone metrics on the complement of the link\,\cite{WY2}, the relative Turaev-Viro invariants of an ideally triangulated $3$-manifold and the volume of the hyperbolic polyhedral metrics on it\,\cite{Y}, and the discrete Fourier transforms of the Yokota invariants of a planar graph with the volume of the deeply truncated polyhedra with the graph as the $1$-skeleton\,\cite{BY}.

On the other hand, the Asymptotic Expansion Conjecture\,\cite{W, J}  relates the perturbative expansion of the Witten-Reshetikhin-Turaev invariants of a closed oriented $3$-manifold to classical invariants of the manifold, including the twisted Reidemeister torsions and the Chern-Simons invariants 
 Recently, it is conjectured  in \cite{O2} and \cite{GRY} that for a closed oriented hyperbolic $3$-manifold $M,$ from the asymptotic expansion of the sequence of the Reshetikhin-Turaev invariants $\{\mathrm {RT}_r(M)\}$ of $M$ evaluated at the root of unity $q=e^{\frac{2\pi\sqrt{-1}}{r}}$ as $r\to \infty$ varying over all the positive odd integers, in addition to the hyperbolic volume and the Chern-Simons invariant of $M,$ 
one should also see the Reidemeister torsion of $M$ twisted by the adjoint action of the holonomy representation of the hyperbolic structure on $M.$
This extends the volume conjecture of the Reshetikhin-Turaev and the Turaev-Viro invariants\,\cite{CY}. See also \cite{OT1, OT2}. 

In this paper, we propose the following Conjecture \ref{CRRT}, Conjecture \ref{CRTV} and Conjecture \ref{CDFT} respectively on the asymptotic expansion of the relative Reshetikhin-Turaev invariants, of the relative Turaev-Viro invariants and of the discrete Fourier transforms of the quantum $6j$-symbols, extending the corresponding volume conjectures proposed in \cite{WY2, BY, Y}. Using the asymptotic estimates of the relative quantum invariants established in \cite{WY2,Y,BY}, in Theorem \ref{main1}, Theorem \ref{main2} and Theorem \ref{main3}, we prove our conjectures for families of cases by reinterpreting the leading terms in the asymptotic expansions  as   geometric quantities associated with  the relevant hyperbolic cone, hyperbolic polyhedral, or deeply truncated tetrahedral metrics of the underlying manifolds or tetrahedra. 
The significance of these expansions is that we do not specify the way that the sequence of colorings converges to the limit. As a consequence, the terms in the expansion will have to depend on $r,$ but the dependence is in a way that  the terms  are purely geometric invariants of the metrics on the underlying manifold and only the metrics vary with $r.$ Besides, new factors involving logarithmic holonomies of longitudes and edge-length terms in respectively the hyperbolic cone metric setting and the polyhedral metric setting naturally appear in the asymptotics expansions of the relative invariants. These factors have not been seen  previously in the asymptotic expansions of the colored Jones polynomials of links or the Reshetikhin-Turaev invariants of closed $3$-manifolds\,\cite{G,W,J,O2,GRY}.

\subsection{Asymptotic expansion of the relative Reshetikhin-Turaev invariants}
Let $M$ be a closed oriented $3$-manifold and let $L$ be a framed hyperbolic link in $M$ with $|L|$ components. Let $\{\mathbf a^{(r)}\}=\{(a^{(r)}_1,\dots,a^{(r)}_{|L|})\}$ be a sequence of colorings of the components of $L$ by the elements of $\{0,\dots,r-2\}$ such that for each $k\in\{1,\dots,|L|\},$ either $a_k^{(r)}> \frac{r}{2}$ for  all $r$ sufficiently large or $a_k^{(r)}< \frac{r}{2}$ for all $r$ sufficiently large. In the former case we let $\mu_k=1$ and in the latter case we let $\mu_k=-1,$ and we let 
$$\theta^{(r)}_k=\mu_k\bigg(\frac{4\pi a^{(r)}_k}{r}-2\pi\bigg).$$
Let $\theta^{(r)}=(\theta^{(r)}_1,\dots,\theta^{(r)}_{|L|}).$ Suppose for all $r$ sufficiently large, a hyperbolic cone metric on $M$ with singular locus $L$ and cone angles $\theta^{(r)}$ exists. We denote $M$ with such a hyperbolic cone metric by $M^{(r)},$  let $\mathrm{Vol}(M^{(r)})$ and $\mathrm{CS}(M^{(r)})$ respectively be the volume and the Chern-Simons invariant of $M^{(r)}$. For $k=1,\dots, |L|$, we arbitrarily choose an orientation of meridian $u_k$ and an orientation of longitude $v_k$ such that their algebraic intersection number $u_k\cdot v_k = 1$, and define the logarithmic holonomy of $u_k$ to be $\mathrm{H}^{(r)}(u_k)=\sqrt{-1}\theta_k^{(r)}$; and for the parallel copies $(\gamma_1,\dots,\gamma_{|L|})$ of the core curves of $L$ given by the framing, let $\mathrm H^{(r)}(\gamma_k)$ be the induced logarithmic holonomy in $M^{(r)}$. 
Let $\rho_{M^{(r)}}:\pi_1(M\setminus L)\to \mathrm{PSL}(2;\mathbb C)$ be the holonomy representation of  the restriction of $M^{(r)}$ to $M\setminus L,$ and let $\mathbb{T}_{(M\setminus L,\mathbf m)}([\rho_{M^{(r)}}])$ be the Reideimester torsion of $M\setminus L$  twisted by the adjoint action of  $\rho_{M^{(r)}}$ with respect to the system of meridians $\mathbf m$ of a tubular neighborhood of the core curves of $L.$

\begin{conjecture}\label{CRRT} Suppose $\{\theta^{(r)}\}$ converges as $r$ tends to infinity.
Then as $r$ varies over all positive odd integers and at $q=e^{\frac{2\pi \sqrt{-1}}{r}},$ the relative Reshetikhin-Turaev invariants 
$$ \mathrm{RT}_r(M,L,\mathbf a^{(r)})=C\frac{e^{\frac{1}{2}\sum_{k=1}^{|L|}\mu_k\mathrm H^{(r)}(\gamma_k)}}{\sqrt{\pm\mathbb{T}_{(M\setminus L,\mathbf m)}([\rho_{M^{(r)}}])}}e^{\frac{r}{4\pi}\big(\mathrm{Vol}(M^{(r)})+\sqrt{-1}\mathrm{CS}(M^{(r)})\big)}\bigg(1+O\Big(\frac{1}{r}\Big)\bigg),$$
where $C$ is a quantity of norm $1$ independent of the geometric structure on $M.$
\end{conjecture}

\begin{theorem} \label{main1}
Conjecture \ref{CRRT} is true if $(M,L)$ is obtained from a fundamental shadow link complement by doing a change-of-pair operation and the limiting cone angles $\theta_1,\dots,\theta_{|L|}$ are sufficiently small.
\end{theorem}


\subsection{Asymptotic expansion of the relative Turaev-Viro invariants}
Suppose $N$ is a $3$-manifold with non-empty boundary and $\mathcal T$ is an ideal triangulation of $N$ with the set of edges $E.$  Let $\{\mathbf b^{(r)}\}=\{(b_1^{(r)},\dots,b_{|E|}^{(r)})\}$ be a sequence of colorings of $(N,\mathcal T)$  by the elements of $\{0,\dots, r-2\}$ such that for each $k\in \{1,\dots, |E|\},$ either $b_k^{(r)} > \frac{r}{2}$ for all $r$ sufficiently large or  $b_k^{(r)} < \frac{r}{2}$ for all $r$ sufficiently large. In the former case we let $\mu_k=1$ and in the latter case we let $\mu_k=-1,$ and we let 
$$\theta^{(r)}_k=\mu_k\bigg( \frac{4\pi b_k^{(r)}}{r}-2\pi\bigg).$$
Let $\theta^{(r)}=(\theta^{(r)}_1,\dots,\theta^{(r)}_{|E|}).$ Suppose for all $r$ sufficiently large, a hyperbolic polyhedral metric on $N$ with cone angles $\theta^{(r)}$ exists. We denote $N$ with such hyperbolic polyhedral metric by $N^{(r)},$  let $\mathrm{Vol}(N^{(r)})$ be the volume of $N^{(r)},$ and let $l^{(r)}_1,\dots, l^{(r)}_{|E|}$ be the lengths of the edges in $N^{(r)}.$ Let $M$ be the $3$-manifold with toroidal boundary obtained from the double of $N$ by removing the double of all the edges, let $\rho_{M^{(r)}}:\pi_1(M)\to \mathrm{PSL}(2;\mathbb C)$ be the holonomy representation of   the restriction of the double of the hyperbolic polyhedral metric on $N^{(r)}$ to $M$ and let $\mathbb{T}_{(M,\mathbf m)}([\rho_{M^{(r)}}])$ be the Reideimester torsion of $M$  twisted by the adjoint action of  $\rho_{M^{(r)}}$ with respect to the system of meridians $\mathbf m$ of a tubular neighborhood of the double of the edges.  

\begin{conjecture}\label{CRTV}
Suppose $\{\theta^{(r)}\}$ converges as $r$ tends to infinity. Then as $r$ varies over all positive odd integers and at $q=e^{\frac{2\pi\sqrt{-1}}{r}},$ the relative Turaev-Viro invariants 
\begin{equation*}
\begin{split}
\mathrm {TV}_r&(N,E,\mathbf b^{(r)})=C \frac{e^{-\sum_{k=1}^{|E|}\mu_kl^{(r)}_k}}{\sqrt{\pm\mathbb{T}_{(M,\mathbf m)}([\rho_{M^{(r)}}])}}  r^{\frac{3}{2}\chi(N)} e^{\frac{r}{2\pi}\mathrm{Vol}(N^{(r)})} \bigg( 1 + O\Big(\frac{1}{r}\Big)\bigg),
\end{split}
\end{equation*}
where $C=\frac{(-1)^{|E|+\chi(N)(\frac{r}{2}-\frac{1}{4})}2^{\mathrm{rank H}_2(N;\mathbb Z_2)}}{(4\pi)^{\chi(N)}}$ 
is a quantity independent of the geometric  structure on $N,$ and $\chi(N)$ is the Euler characteristic of $N.$
\end{conjecture}

\begin{theorem} \label{main2}  
Conjecture \ref{CRTV} is true if the limiting cone angles $ \theta_1,\dots,\theta_{|E|}$ are sufficiently small. 
\end{theorem}


\subsection{Asymptotic expansion of the discrete Fourier transforms of quantum $6j$-symbols}

Let $(I,J)$ be a partition of $\{1,\dots,6\},$ and let $\{(\mathbf b_I^{(r)}, \mathbf a_J^{(r)})\}=\{((b_i^{(r)})_{i\in I}, (a_j^{(r)})_{i\in J})\}$ be a sequence of $6$-tuples of the elements of $\{0,\dots, r-2\}$  such that for any $i\in I,$ either $b_i^{(r)}>\frac{r}{2}$ for all $r$ or $b_i^{(r)}<\frac{r}{2}$ for all $r;$ and for any $j\in J,$ either $a_j^{(r)}>\frac{r}{2}$ for all $r$ or $a_j^{(r)}<\frac{r}{2}$ for all $r.$ In the former case we let $\mu_i=\mu_j=1$ and in the latter case we let $\mu_i=\mu_j=-1,$ and for $i\in I$ we let
$$\theta^{(r)}_i=\mu_i\bigg(\frac{2\pi b_i^{(r)}}{r}-\pi\bigg)$$
and for $j\in J$ we let
$$\theta^{(r)}_j=\mu_j\bigg(\frac{2\pi a_j^{(r)}}{r}-\pi\bigg).$$ 
For each $r,$ suppose $\Delta^{(r)}$ is the deeply truncated tetrahedron  of type $(I,J)$ with  the set of dihedral angles at the edges of deep truncation $\theta^{(r)}_I=(\theta^{(r)}_i)_{i\in I}$ and the set of dihedral angles at the regular edges $\theta^{(r)}_J=(\theta^{(r)}_j)_{j\in J}.$ We let $\mathrm{Vol}(\Delta^{(r)})$ be the volume of $\Delta^{(r)}.$
Let $l^{(r)}_I=(l^{(r)}_i)_{i\in I}$ be the set of the lengths of the edges of deep truncation and let $(l^{(r)}_j)_{j\in J}$ be the set of the lengths of the regular edges, and let $\mathbb G\big(l^{(r)}_I,\sqrt{-1}\theta^{(r)}_J\big)$ be the value of the Gram matrix function at $(l^{(r)}_I,\sqrt{-1}\theta^{(r)}_J).$  

\begin{conjecture}\label{CDFT}
Suppose $\{\theta^{(r)}_I\}$ and $\{\theta^{(r)}_J\}$ converge as $r$ tends to infinity. Then as $r$ varies over all positive odd integers and evaluated at the root of unity $q=e^{\frac{2\pi \sqrt{-1}}{r}},$ the discrete Fourier transforms of the Yokota invariant of the trivalent graph $\cp{\includegraphics[width=0.5cm]{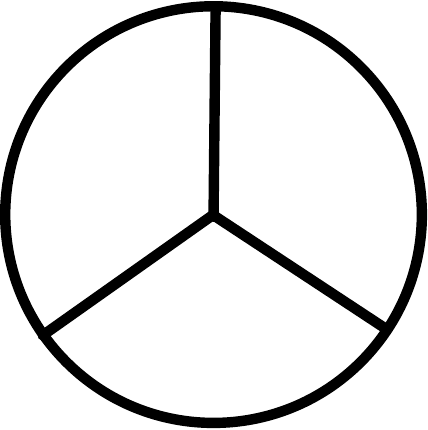}}$ 
\begin{equation*}
\mathrm{\widehat {Y}}_r(\mathbf b^{(r)}_I; \mathbf a^{(r)}_J)=C\frac{e^{-\sum_{k=1}^{6}\mu_kl^{(r)}_k}}{\sqrt{-\det\bigg(\frac{\partial \theta_{i_1}}{\partial l_{i_2}}\Big|_{\big(l^{(r)}_I,\theta^{(r)}_J\big)}\bigg)_{i_1,i_2\in I}\det\mathbb G\big(l^{(r)}_I,\sqrt{-1}\theta^{(r)}_J\big)}}  r^{\frac{3 |I|-6}{2}}e^{\frac{r}{\pi}\mathrm{Vol}(\Delta^{(r)})} \bigg( 1 + O\Big(\frac{1}{r}\Big)\bigg),
\end{equation*}
where $C=\frac{(-1)^{\frac{3}{2}+\frac{r}{2}(|I|-2)}n(a_J)}{2^{\frac{3|I|}{2}-1}\pi^{|I|-2}}$ is a quantity independent of the geometry of the deeply truncated tetrahedron, and $n(a_J)$ is the number of $3$-admissible colorings that have the same parity as $\mathbf a_J.$
\end{conjecture}

\begin{theorem} \label{main3}  
Conjecture \ref{CDFT} is true if the limiting dihedral angles $\theta_1,\dots,\theta_6$ are sufficiently small. 
\end{theorem}

\begin{remark}  In the case that $I=\emptyset,$ one should compare Theorem \ref{main3} with Chen-Murakami\,\cite[Theorem 2]{CM}, where they consider the asymptotic expansion of the norm of the quantum $6j$-symbols with a specific choice of the  sequence of the colorings.
\end{remark}
\bigskip

\noindent\textbf{Acknowledgments.}  The authors would like to thank Qingtao Chen for explaining his joint work with Jun Murakami. The authors are also grateful to the referee for reading the paper carefully and providing valuable suggestions. The second author is partially supported by NSF Grant DMS-1812008.


\section{Preliminaries}

\subsection{Adjoint Twisted Reidemeister torsion and Gram matrices}

Let $\mathrm C_*$ be a finite chain complex 
$$0\to \mathrm C_d\xrightarrow{\partial}\mathrm C_{d-1}\xrightarrow{\partial}\cdots\xrightarrow{\partial}\mathrm C_1\xrightarrow{\partial} \mathrm C_0\to 0$$
of $\mathbb C$-vector spaces, and for each $\mathrm C_k$ choose a basis $\mathbf c_k.$ Let $\mathrm H_*$ be the homology of $\mathrm C_*,$ and for each $\mathrm H_k$ choose a basis $\mathbf h_k$ and a lift $\widetilde{\mathbf h}_k\subset \mathrm C_k$ of $\mathbf h_k.$ We also choose a basis $\mathbf b_k$ for each image $\partial (\mathrm C_{k+1})$ and a lift $\widetilde{\mathbf b}_k\subset \mathrm C_{k+1}$ of $\mathbf b_k.$ Then $\mathbf b_k\sqcup \widetilde{\mathbf b}_{k-1}\sqcup \widetilde{\mathbf h}_k$ form a basis of $\mathrm C_k.$ Let $[\mathbf b_k\sqcup \widetilde{\mathbf b}_{k-1}\sqcup \widetilde{\mathbf h}_k;\mathbf c_k]$ be the determinant of the transition matrix from the standard basis $\mathbf c_k$ to the new basis $\mathbf b_k\sqcup \widetilde{\mathbf b}_{k-1}\sqcup \widetilde{\mathbf h}_k.$
 Then the Reidemeister torsion of the chain complex $\mathrm C_*$ with the chosen bases $\mathbf c_*$ and $\mathbf h_*$ is defined by 
\begin{equation*}
\mathrm{Tor}(\mathrm C_*, \{\mathbf c_k\}, \{\mathbf h_k\})=\pm\prod_{k=0}^d[\mathbf b_k\sqcup \widetilde{\mathbf b}_{k-1}\sqcup \widetilde{\mathbf h}_k;\mathrm c_k]^{(-1)^{k+1}}.
\end{equation*}
It is easy to check that $\mathrm{Tor}(\mathrm C_*, \{\mathbf c_k\}, \{\mathbf h_k\})$ depends only on the choice of $\{\mathbf c_k\}$ and $\{\mathbf h_k\},$ and does not depend on the choices of $\{\mathbf b_k\}$ and the lifts $ \{\widetilde{\mathbf b}_k\}$ and $\{\widetilde{\mathbf h}_k\}.$

We recall the twisted Reidemeister torsion of a CW-complex following the conventions in \cite{P2}. Let $K$ be a finite CW-complex and let $\rho:\pi_1(M)\to\mathrm{SL}(N;\mathbb C)$ be a representation of its fundamental group. Consider the twisted chain complex 
$$\mathrm C_*(K;\rho)= \mathbb C^N\otimes_\rho \mathrm C_*(\widetilde K;\mathbb Z)$$
where $\mathrm C_*(\widetilde K;\mathbb Z)$ is the simplicial complex of the universal covering of $K$ and $\otimes_\rho$ means the tensor product over $\mathbb Z$ modulo the relation
$$\mathbf v\otimes( \gamma\cdot\mathbf c)=\Big(\rho(\gamma)^T\cdot\mathbf v\Big)\otimes \mathbf c,$$
where $T$ is the transpose, $\mathbf v\in\mathbb C^N,$ $\gamma\in\pi_1(K)$ and $\mathbf c\in\mathrm C_*(\widetilde K;\mathbb Z).$ The boundary operator on $\mathrm C_*(K;\rho)$ is defined by
$$\partial(\mathbf v\otimes \mathbf c)=\mathbf v\otimes \partial(\mathbf c)$$
for $\mathbf v\in\mathbb C^N$ and $\mathbf c\in\mathrm C_*(\widetilde K;\mathbb Z).$ Let $\{\mathbf e_1,\dots,\mathbf e_N\}$ be the standard basis of $\mathbb C^N,$ and let $\{c_1^k,\dots,c_{d^k}^k\}$ denote the set of $k$-cells of $K.$ Then we call
$$\mathbf c_k=\big\{ \mathbf e_i\otimes c_j^k\ \big|\ i\in\{1,\dots,N\}, j\in\{1,\dots,d^k\}\big\}$$
the standard basis of $\mathrm C_k(K;\rho).$ Let $\mathrm H_*(K;\rho)$ be the homology of the chain complex $\mathrm C_*(K;\rho)$ and let $\mathbf h_k$ be a basis of $\mathrm H_k(K;\rho).$ Then the Reidemeister torsion of $K$ twisted by $\rho$ with basis $\{\mathbf h_k\}$ is 
$$\mathrm{Tor}(K, \{\mathbf h_k\}; \rho)=\mathrm{Tor}(\mathrm C_*(K;\rho),\{\mathbf c_k\}, \{\mathbf h_k\}).$$

By \cite{P}, $\mathrm{Tor}(K, \{\mathbf h_k\}; \rho)$ depends only on the conjugacy class of $\rho.$ By for e.g. \cite{T2}, the Reidemeister torsion is invariant under elementary expansions and elementary collapses of CW-complexes, and by \cite{M}  it is invariant under subdivisions, hence defines an invariant of PL-manifolds and of topological manifolds of dimension less than or equal to $3.$

We list some results by Porti\,\cite{P} for the Reidemeister torsions of hyperbolic $3$-manifolds twisted by the adjoint representation $\mathrm {Ad}_\rho=\mathrm {Ad}\circ\rho$ of the holonomy  $\rho$ of the hyperbolic structure. Here $\mathrm {Ad}$ is the adjoint acton of $\mathrm {PSL}(2;\mathbb C)$ on its Lie algebra $\mathbf{sl}(2;\mathbb C)\cong \mathbb C^3.$

For a closed oriented hyperbolic $3$-manifold  $M$ with the holonomy representation $\rho,$ by the Weil local rigidity theorem and the Mostow rigidity theorem,
$\mathrm H_k(M;\mathrm{Ad}_\rho)=0$ for all $k.$ Then the twisted Reidemeister torsion 
$$\mathrm{Tor}(M;\mathrm{Ad}_\rho)\in\mathbb C^*/\{\pm 1\}$$
 is defined without making any additional choice.

For a compact, orientable  $3$-manifold  $M$ with boundary consisting of $n$ disjoint tori $T_1 \dots,  T_n$ whose interior admits a complete hyperbolic structure with  finite volume, let $\mathrm X(M)$ be the $\mathrm{SL}(2;\mathbb C)$-character variety of $M,$ let $\mathrm X_0(M)\subset\mathrm X(M)$ be the distinguished component containing the character of a chosen lift of the holonomy representation of the complete hyperbolic structure of $M,$ and let $\mathrm X^{\text{irr}}(M)\subset\mathrm X(M)$ be the subset consisting of the irreducible characters. 
 
\begin{theorem}\cite[Section 3.3.3]{P}\label{HM} For a generic character $[\rho]\in\mathrm X_0(M)\cap\mathrm X^{\text{irr}}(M)$  we have:
\begin{enumerate}[(1)]
\item For $k\neq 1,2,$ $\mathrm H_k(M;\mathrm{Ad}\rho)=0.$
\item  For $i\in\{1,\dots,n\},$ let $\mathbf I_i\in \mathbb C^3$ be up to scalar the unique invariant vector of $\mathrm Ad_\rho(\pi_1(T_i)).$ Then
$$\mathrm H_1(M;\mathrm{Ad}\rho)\cong\bigoplus_{i=1}^n\mathrm H_1(T_i;\mathrm{Ad}\rho)\cong \mathbb C^n.$$ 
Moreover, for each $\alpha=([\alpha_1],\dots,[\alpha_n])\in \mathrm H_1(\partial M;\mathbb Z)\cong 
\bigoplus_{i=1}^n\mathrm H_1(T_i;\mathbb Z)$ with $[\alpha_1],\dots,[\alpha_n] \neq 0$, $\mathrm H_1(M;\mathrm{Ad}\rho)$  has a basis 
$$\mathbf h^1_{(M,\alpha)}=\{\mathbf I_1\otimes [\alpha_1],\dots, \mathbf I_n\otimes [\alpha_n]\}.$$
\item Let $([T_1],\dots,[T_n])\in \bigoplus_{i=1}^n\mathrm H_2(T_i;\mathbb Z)$ be the fundamental classes of $T_1,\dots, T_n.$ Then 
 $$\mathrm H_2(M;\mathrm{Ad}\rho)\cong\bigoplus_{i=1}^n\mathrm H_2(T_i;\mathrm{Ad}\rho)\cong \mathbb C^n,$$ 
and has a basis 
$$\mathbf h^2_M=\{\mathbf I_1\otimes [T_1],\dots, \mathbf I_n\otimes [T_n]\}.$$
\end{enumerate}
\end{theorem}

\begin{remark}[\cite{P}] Important examples of the generic characters in Theorem \ref{HM} include the characters of the lifting in $\mathrm{SL}(2;\mathbb C)$ of the holonomy of the complete hyperbolic structure on the interior of $M,$
 the restriction of the holonomy of the closed $3$-manifold $M_\mu$ obtained from $M$ by doing the hyperbolic Dehn surgery along the system of simple closed curves $\mu$ on $\partial M,$
 and by \cite{HK} the holonomy of a hyperbolic structure on the interior of $M$ whose completion is a conical manifold with cone angles less than $2\pi.$
\end{remark}

For $\alpha\in\mathrm H_1(M;\mathbb Z),$ define $\mathbb T_{(M,\alpha)}$ on $\mathrm X_0(M)$ by
$$\mathbb T_{(M,\alpha)}([\rho])=\mathrm{Tor}(M, \{\mathbf h^1_{(M,\alpha)},\mathbf h^2_M\};\mathrm{Ad}_\rho)$$
for any generic $[\rho]\in \mathrm X_0(M)\cap\mathrm X^{\text{irr}}(M)$ in Theorem \ref{HM}, and equals $0$ otherwise.

\begin{theorem}\cite[Theorem 4.1]{P}\label{funT}
Let $M$ be a compact, orientable  $3$-manifold with boundary consisting of $n$ disjoint tori $T_1 \dots,  T_n$ whose interior admits a complete hyperbolic structure with  finite volume. Let  $\mathbb C(\mathrm X_0(M))$ be the ring of rational functions over $\mathrm X_0(M).$ Then there is up to sign a unique function
\begin{equation*}
\begin{split}
\mathrm H_1(\partial M;\mathbb Z)&\to \mathbb C(\mathrm X_0(M))\\
\alpha\quad\quad &\mapsto \quad\mathbb T_{(M,\alpha)}
\end{split}
\end{equation*}
which is a $\mathbb Z$-multilinear homomorphism with respect to the direct sum $\mathrm H_1(\partial M;\mathbb Z)\cong 
\bigoplus_{i=1}^n\mathrm H_1(T_i;\mathbb Z)$ satisfying the following properties:
\begin{enumerate}[(i)]
\item For all $\alpha \in \mathrm H_1(\partial M;\mathbb Z),$ the domain of definition of $\mathbb T_{(M,\alpha)}$ contains an open subset $\mathrm X_0(M)\cap\mathrm X^{\text{irr}}(M).$

\item \emph{(Change of curves formula).} Let $\mu=\{\mu_1,\dots,\mu_n\}$ and $\gamma=\{\gamma_1,\dots,\gamma_n\}$ be two systems of simple closed curves on $\partial M.$ If $\mathrm H(\mu_1),\dots, \mathrm H(\mu_n)$ and $\mathrm H(\gamma_1),\dots,\mathrm H(\gamma_n)$ are respectively the logarithmic holonomies of the curves in $\mu$ and $\gamma,$ then we have the equality of rational functions
\begin{equation*}\label{coc}
\mathbb T_{(M,\mu)}
=\pm\det\bigg( \frac{\partial \mathrm H(\mu_i)}{\partial \mathrm H(\gamma_j)}\bigg)_{ij}\mathbb T_{(M,\gamma)}.
\end{equation*}
\item \emph{(Surgery formula).} Let $[\rho_\mu]\in \mathrm X_0(M)$ be the character induced by the holonomy of the closed $3$-manifold $M_\mu$ obtained from $M$ by doing the hyperbolic Dehn surgery along the system of simple closed curves $\mu$ on $\partial M.$ If $\mathrm H(\gamma_1),\dots,\mathrm H(\gamma_n)$ are the logarithmic holonomies of the core curves $\gamma_1,\dots,\gamma_n$ of the solid tori added. Then
\begin{equation*}\label{sf}
\mathrm{Tor}(M_\mu;\mathrm{Ad}_{\rho_\mu})=\pm\mathbb T_{(M,\mu)}([\rho_\mu])\prod_{i=1}^n\frac{1}{4\sinh^2\frac{\mathrm H(\gamma_i)}{2}}.
\end{equation*}
\end{enumerate}
\end{theorem}

Next we list some results for the computation of twisted Reidemeister torsions from \cite{WY3}. We first recall that if $\mathrm M_{4\times 4}(\mathbb C)$ is the space of $4\times 4$ matrices with complex entries, then the \emph{Gram matrix function} 
$$\mathbb G:\mathbb C^6\to \mathrm M_{4\times 4}(\mathbb C)$$
 is defined by 
\begin{equation}\label{gram}
\begin{split}
\mathbb{G}(\mathbf z)=\left[
\begin{array}{cccc}
1 & -\cosh z_{1} & -\cosh z_{2} &-\cosh z_{6}\\
-\cosh z_{1}& 1 &-\cosh z_{3} & -\cosh z_{5}\\
-\cosh z_{2} & -\cosh z_{3} & 1 & -\cosh z_{4} \\
-\cosh z_{6} & -\cosh z_{5} & -\cosh z_{4}  & 1 \\
 \end{array}\right]
 \end{split}
\end{equation}
for $\mathbf z=(z_{1}, z_{2}, z_{3}, z_{4}, z_{5}, z_{6})\in\mathbb C^6.$ The values of $\mathbb G$ at different $\mathbf z$ recover the Gram matrices of deeply truncated tetrahedra  of all the types. See Section \ref{dtt} or \cite[Section 2.1]{BY} for more details.

\begin{theorem}\cite[Theorem 1.1]{WY3} \label{Tor1} Let $M=\#^{c+1}(S^2\times S^1)\setminus L_{\text{FSL}}$ be the complement of a fundamental shadow link $L_{\text{FSL}}$ with $n$ components $L_1,\dots,L_n,$ which is the orientable double of the union of truncated tetrahedra $\Delta_1,\dots, \Delta_c$ along pairs of the triangles of truncation, and let $\mathrm X_0(M)$ be the distinguished component of the $\mathrm{SL}(2;\mathbb C)$ character variety of $M$ containing a lifting of the holonomy representation of the complete hyperbolic structure.  
\begin{enumerate}[(1)] 
\item  Let $\mathbf u=(u_1,\dots, u_n)$ be the system of the meridians of a tubular neighborhood of the components of $L_{\text{FSL}}.$ For a generic irreducible character $[\rho]$ in $\mathrm X_0(M),$ let $\mathrm H(u_1),
\dots,\mathrm H(u_n)$ be the logarithmic holonomies of $\mathbf u.$ For each $s\in\{1,\dots,c\},$ let $L_{s_1},\dots,L_{s_6}$ be the components of $L_{\text{FSL}}$ intersecting $\Delta_s,$ and let $\mathbb G_s$ be the value of the Gram matrix function at $\Big(\frac{\mathrm H(u_{s_1})}{2},\dots,\frac{\mathrm H(u_{s_6})}{2}\Big).$ Then
 $$\mathbb T_{(M,\mathbf u)}([\rho])=\pm2^{3c}\prod_{s=1}^c \sqrt{\det\mathbb G_s}.$$

\item In addition to the assumptions and notations of (1), let $\mathbf m=(m_1,\dots,m_n)$ be a system of simple closed curves on $\partial M,$ and let $(\mathrm H(m_1),\dots, \mathrm H(m_n))$ be their logarithmic holonomies which are functions of $(\mathrm H(u_1),
\dots,\mathrm H(u_n)).$ Then
 $$ \mathbb T_{(M,\mathbf m)}[(\rho)]=\pm2^{3c}\det\bigg(\frac{\partial \mathrm H(m_i)}{\partial \mathrm H(u_j)}\bigg|_{[\rho]}\bigg)_{ij}\prod_{s=1}^c \sqrt{\det\mathbb G_s}.$$

\item Suppose $M_{\mathbf m}$ is the closed $3$-manifold obtained from $M$ by doing the hyperbolic Dehn surgery along a system of simple closed curves $\mathbf m=(m_1,\dots,m_n)$ on $\partial M$ and $\rho_{\mathbf m}$ is the restriction of the holonomy representation of $M_{\mathbf m}$ to $M.$ Let $(\mathrm H(m_1),\dots, \mathrm H(m_n))$ be the logarithmic holonomies of $\mathbf m$ which are functions of the logarithmic holonomies of the meridians $\mathbf u.$ Let $(\gamma_1,\dots,\gamma_n)$ be a system of simple closed curves on $\partial M$ that are isotopic to the core curves of the solid tori filled in and let $\mathrm H(\gamma_1),\dots,\mathrm H(\gamma_n)$ be their logarithmic holonomies in $[\rho_\mu].$ Let  $\mathrm H(u_1),
\dots,\mathrm H(u_n)$ be the logarithmic holonomies of the meridians $\mathbf u$ in $[\rho_\mu]$ and for each $s\in\{1,\dots,c\},$ let $L_{s_1},\dots,L_{s_6}$ be the components of $L_{\text{FSL}}$ intersection $\Delta_s$ and let $\mathbb G_s$ be the value of the Gram matrix function at $\Big(\frac{\mathrm H(u_{s_1})}{2},\dots,\frac{\mathrm H(u_{s_6})}{2}\Big).$
Then 
 $$ \mathrm{Tor}(M_{\mathbf m};\mathrm{Ad}_{\rho_{\mathbf m}})=\pm2^{3c-2n}\det\bigg(\frac{\partial \mathrm H(m_i)}{\partial \mathrm H(u_j)}\bigg|_{[\rho_\mu]}\bigg)_{ij}\prod_{s=1}^c \sqrt{\det\mathbb G_s}\prod_{i=1}^n\frac{1}{\sinh^2\frac{\mathrm H(\gamma_i)}{2}}.$$
\end{enumerate} 
\end{theorem}

\begin{theorem}\cite[Theorem 1.4]{WY3} \label{Tor2} Let $M$ be the double of a hyperbolic polyhedral $3$-manifold $N$ which is the union of truncated tetrahedra $\Delta_1,\dots,\Delta_{|T|}$ along pairs of hexagonal faces with the double of the edges $e_1,\dots,e_{|E|}$ removed. For $i\in \{1,\dots,{|E|}\},$ let $l_i$ be the lengths of $e_i.$ 
\begin{enumerate}[(1)] 
\item  Let $\mathbf u$ be the system of the preferred longitudes of $\partial M$ with the logarithmic holonomies $2l_1,$ $\dots,$ $2l_{|E|}.$ For each $s\in\{1,\dots,{|T|}\},$ let $e_{s_1},\dots,e_{s_6}$ be the edges intersecting $\Delta_s,$ and let $\mathbb G_s$ be the value of the Gram matrix function at $( l_{s_1},\dots, l_{s_6}).$ Let $\rho$ be the holonomy representation of the hyperbolic cone metric on $M$ obtained from the double of the hyperbolic polyhedral metric of $N.$
Then
 $$\mathbb T_{(M,\mathbf u)}([\rho])=\pm2^{3|T|}\prod_{s=1}^{|T|} \sqrt{\det\mathbb G_s}.$$

\item In addition to the assumptions and notations of (1), let $\mathbf m$ be the system of the meridians of tubular neighborhoods of the double of the edges, and let $(\theta_1,\dots,\theta_n)$ be the cone angle functions of the lengths of $N.$ Then
$$ \mathbb T_{(M,\mathbf m)}([\rho])=\pm(-1)^{\frac{3|E|}{2}}2^{3|T|-|E|}\det\bigg(\frac{\partial \theta_i}{\partial l_j}\bigg|_{[\rho]}\bigg)_{ij}\prod_{s=1}^{|T|} \sqrt{\det\mathbb G_s}.$$

\item Suppose $\overline {M}$ is the double of a geometrically triangulated hyperbolic $3$-manifold $N$ with totally geodesic boundary (which is $M$ with the removed double of edges filled back) and $\overline\rho$ is its holonomy representation.
Let $(\theta_1,\dots,\theta_n)$ be the cone angle functions in terms of the edge lengths of $N,$ and let  $l_1,\dots,l_n$ be the lengths of the edges of $N$ in $\overline \rho.$  For each $s\in\{1,\dots,d\},$ let $e_{s_1},\dots,e_{s_6}$ be the edges intersecting $\Delta_s$ and let $\mathbb G_s=\mathbb G( l_{s_1},\dots, l_{s_6})$ be the value of the Gram matrix function at $( l_{s_1},\dots, l_{s_6}).$ Then 
 $$ \mathrm{Tor}(\overline {M};\mathrm{Ad}_{\overline\rho})=\pm(-1)^{\frac{3|E|}{2}}2^{3|T|-3|E|}\det\bigg(\frac{\partial \theta_i}{\partial l_j}\bigg|_{[\overline\rho]}\bigg)_{ij}\prod_{s=1}^{|T|} \sqrt{\det\mathbb G_s}\prod_{i=1}^{|E|}\frac{1}{\sinh^2l_i}.$$
\end{enumerate} 
\end{theorem}


\subsection{Deeply truncated tetrahedra}\label{dtt}

A \emph{ deeply truncated tetrahedron} is a compact hyperbolic polyhedron with faces $H_1,$ $H_2,$ $H_3,$ $H_4,$ $T_1,$ $T_2,$ $T_3,$ and $T_4$ such that: (1) For each  $i\in \{1,2,3,4\},$ $T_i\cap H_i=\emptyset.$ (2) For each $\{i,j\} \subset \{1,2,3,4\},$ $T_i\cap H_j\neq\emptyset,$ and the dihedral angle between them is always $\frac{\pi}{2}.$ (3) For each $\{i,j\} \subset \{1,2,3,4\},$ either $T_i\cap T_j\neq\emptyset$ or $H_i\cap H_j\neq\emptyset,$ but not both.

From the definition, we see that each face $H_i$ or $T_i$ is one of the following four types: (1) a hyperbolic triangle, (2) a hyperbolic quadrilateral with two right angles, (3) a hyperbolic pentagon with four right angles and (4) a hyperbolic hexagon with six right angles.

\begin{figure}[htbp]
\centering
\includegraphics[scale=0.25]{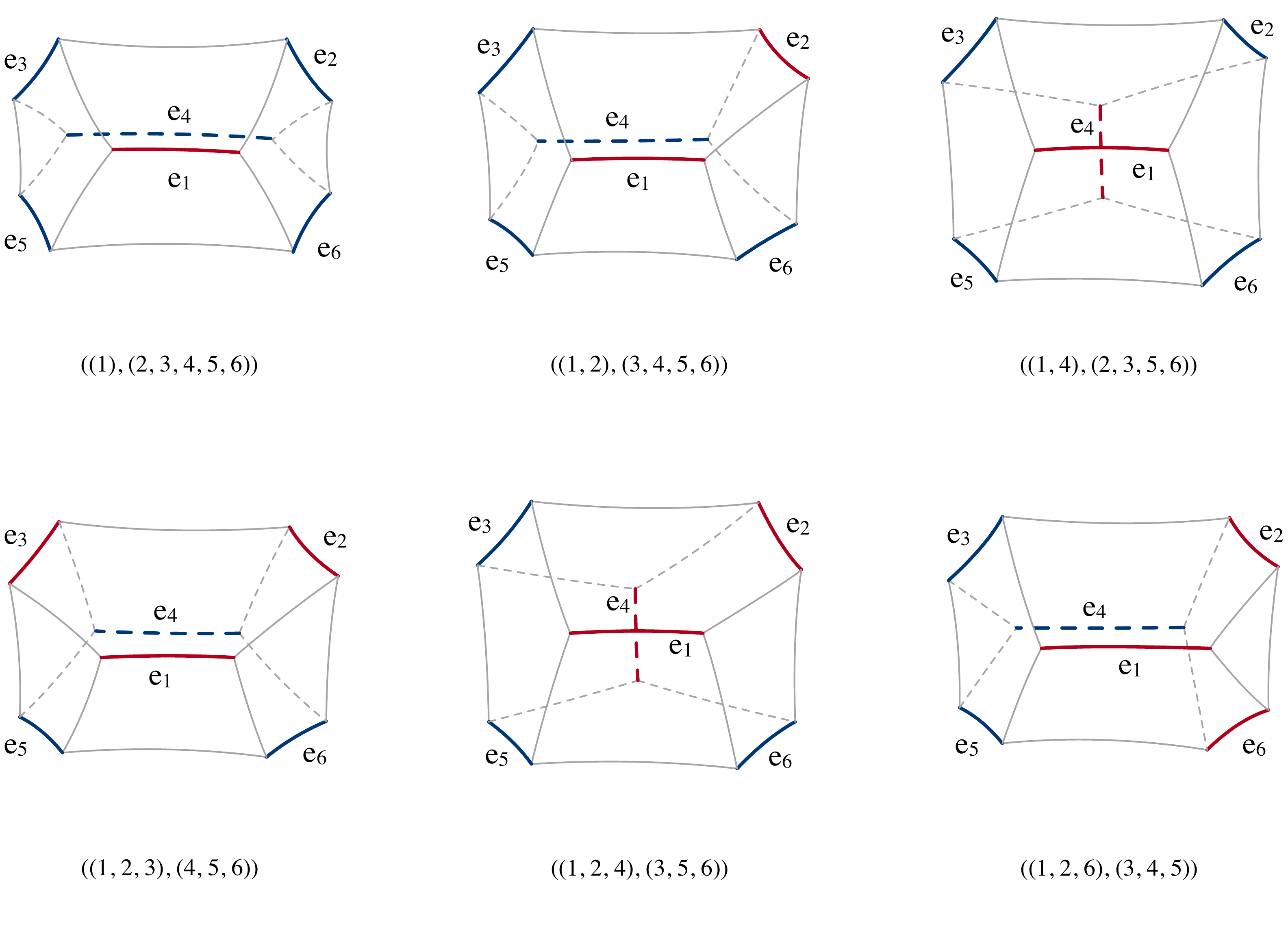}
\caption{Deeply truncated tetrahedra: The red edges are the edges of deep truncation and the blue edges are the regular edges. The dihedral angles at the grey edges are all right angles.}
\label{deep}
\end{figure}

We only consider the intersection of $H_i$ and $H_j$ or the intersection of $T_i$ and $T_j$ as the \emph{edge} of the deeply truncated tetrahedron; therefore there are in total six edges. We call an edge between $H_i$ an $H_j$ a \emph{regular edge} and an edge between $T_i$ and $T_j$ an \emph{edge of deep truncation}. Let $(I,J)$ be a partition of $\{1,\dots,6\}.$ A deeply truncated tetrahedron $\Delta$ is of type $(I,J)$ if $\{e_i\}_{i\in I}$ is the set of edges of deep truncation.  Up to a permutation of the indices, all the types of deeply truncated tetrahedra besides the truncated hyperideal tetrahedron are listed in Figure \ref{deep}.

For a deeply truncated tetrahedron $\Delta$ of type $(I,J),$ let $\{l_i\}_{i\in I}$ be the lengths of the edges of deep truncation and let  $\{\theta_j\}_{j\in J}$ be the
dihedral angles at the regular edges. Let $c_i=\cosh l_i$ for $i\in I$ and $c_j=\cos\theta_j$ for $j\in J.$ Then  the \emph{Gram matrix} of $\Delta$ is the following matrix 
 \begin{equation}\label{gram2}
\left[
\begin{array}{cccc}
1 & -c_1 & -c_2 & -c_6\\
-c_1& 1 & -c_3 & -c_5\\
-c_2 & -c_3 & 1 & -c_4 \\
-c_6 & -c_5 & -c_4  & 1 \\
 \end{array}\right],
\end{equation}
which equals the value $\mathbb G\big((l_i)_{i\in I},(\sqrt{-1}\theta_j)_{j\in J})$ of the Gram matrix function $\mathbb G$ defined in (\ref{gram}) at $\big((l_i)_{i\in I},(\sqrt{-1}\theta_j)_{j\in J}).$


\subsection{Potential functions for hyperideal and deeply truncated tetrahedra}
We recall the definition and properties of the analytic function $U$ defined in \eqref{term} below, which plays an important role in establishing the asymptotics of the relevant quantum invariants \cite{BY, WY2, Y}. We refer the reader to \cite{BY, WY2, Y} for the relationship between $U$ and the asymptotics of quantum 6j-symbols.

\begin{definition}
A triple $(\alpha_1,\alpha_2,\alpha_3)\in [0,2\pi]^3$ is \emph{admissible} if 
\begin{enumerate}[(1)]
\item $\alpha_i+\alpha_j-\alpha_k\geqslant 0$ for $\{i,j,k\}=\{1,2,3\},$
\item $\alpha_1+\alpha_2+\alpha_3\leqslant 4\pi.$
\end{enumerate}
A $6$-tuple $(\alpha_1,\dots,\alpha_6)\in [0,2\pi]^6$ is \emph{admissible} if the triples $\{1,2,3\},$ $\{1,5,6\},$ $\{2,4,6\}$ and $\{3,4,5\}$ are admissible.
\end{definition}

\begin{definition} A $6$-tuple $(\alpha_1,\dots,\alpha_6)\in [0,2\pi]^6$ is of the \emph{hyperideal type} if for $\{i,j,k\}=\{1,2,3\},$ $\{1,5,6\},$ $\{2,4,6\}$ and $\{3,4,5\},$

\begin{enumerate}[(1)]
\item $0\leqslant \alpha_i+\alpha_j-\alpha_k\leqslant 2\pi,$
\item $2\pi\leqslant \alpha_i+\alpha_j+\alpha_k\leqslant 4\pi.$
\end{enumerate}
\end{definition}

Given an admissible 6-tuple $(\alpha_1,\dots,\alpha_6)\in [0,2\pi]^6$, we let $$\tau_1=\frac{\alpha_1+\alpha_2+\alpha_3}{2}, \ \tau_2=\frac{\alpha_1+\alpha_5+\alpha_6}{2}, \ \tau_3=\frac{\alpha_2+\alpha_4+\alpha_6}{2}, \ \tau_4=\frac{\alpha_3+\alpha_4+\alpha_5}{2},$$ $$\eta_1=\frac{\alpha_1+\alpha_2+\alpha_4+\alpha_5}{2},\ 
\eta_2=\frac{\alpha_1+\alpha_3+\alpha_4+\alpha_6}{2}, \ \eta_3=\frac{\alpha_2+\alpha_3+\alpha_5+\alpha_6}{2}.$$
On the region 
\begin{equation*}
\mathrm{B_{H,\mathbb C}}=\bigg\{(\alpha_1,\dots,\alpha_6,\xi)\in\mathbb C^7\left |
\begin{array}{c}
(\mathrm{Re}(\alpha_1),\dots,\mathrm{Re}(\alpha_6))\text{ is of the hyperideal type,}\\
 \max\{\mathrm{Re}(\tau_i)\}\leqslant \mathrm{Re}(\xi)\leqslant \min\{\mathrm{Re}(\eta_j), 2\pi\}\\
 \end{array}\right\},
 \end{equation*}
we define holomoprhic functions $U$ and $\kappa$ by
\begin{equation}\label{term}
\begin{split}
U(\alpha_1,\dots,\alpha_6,\xi)=&\,\pi^2+\frac{1}{2}\sum_{i=1}^4\sum_{j=1}^3(\eta_j-\tau_i)^2-\frac{1}{2}\sum_{i=1}^4(\tau_i-\pi)^2\\
&+(\xi-\pi)^2-\sum_{i=1}^4(\xi-\tau_i)^2-\sum_{j=1}^3(\eta_j-\xi)^2\\
&-2\mathit{Li}_2(1)-\frac{1}{2}\sum_{i=1}^4\sum_{j=1}^3\mathit{Li}_2\big(e^{2i(\eta_j-\tau_i)}\big)+\frac{1}{2}\sum_{i=1}^4\mathit{Li}_2\big(e^{2i(\tau_i-\pi)}\big)\\
&-\mathit{Li}_2\big(e^{2i(\xi-\pi)}\big)+\sum_{i=1}^4\mathit{Li}_2\big(e^{2i(\xi-\tau_i)}\big)+\sum_{j=1}^3\mathit{Li}_2\big(e^{2i(\eta_j-\xi)}\big)\\
\end{split}
\end{equation}
and 
 \begin{equation}\label{kappa}
 \begin{split}
\kappa(\alpha_1,\dots,\alpha_6,\xi)=&\,\frac{1}{2}\sum_{i=1}^4 \sqrt{-1}\tau_{i}- \sqrt{-1}\xi-\frac{3}{2}\sqrt{-1}\pi\\
&+\frac{1}{4}\sum_{i=1}^4\sum_{j=1}^3\log\big(1-e^{2\sqrt{-1}(\eta_{j}-\tau_{i})}\big)-\frac{3}{4}\sum_{i=1}^4\log\big(1-e^{2\sqrt{-1}(\tau_{i}-\pi)}\big)\\
&+\frac{3}{2}\log\big(1-e^{2\sqrt{-1}(\xi-\pi)}\big)-\frac{1}{2}\sum_{i=1}^4\log\big(1-e^{2\sqrt{-1}(\xi-\tau_{i})}\big)\\
&-\frac{1}{2}\sum_{j=1}^3\log\big(1-e^{2\sqrt{-1}(\eta_{j}-\xi)}\big).
 \end{split}
 \end{equation}
The partial derivatives of $U$ will play a key role in the proof of our main results. First, a direct computation shows that 
\begin{align}\label{partialUxi}
    \frac{\partial U}{\partial \xi}
    =&\ -12\xi - 2\pi + 4\sum_{k=1}^6 \alpha_k + 2\sqrt{-1}\log\left(1-e^{2\sqrt{-1}\xi}\right) \notag\\
    &\ -2\sqrt{-1}\sum_{i=1}^4 \log\left(1-e^{2\sqrt{-1}(\xi-\tau_i)}\right) +2\sqrt{-1}
    \sum_{j=1}^3\log\left(1-e^{2\sqrt{-1}(\eta_j-\xi)}\right).
\end{align}
Besides, we have the following 8 equalities.
\begin{align*}
    \sum_{k=1}^6\frac{\partial}{\partial \alpha_k}
    \left(\frac{1}{2}\sum_{i=1}^4\sum_{j=1}^3(\eta_j-\tau_i)^2\right)
    &= \frac{1}{2}\sum_{k=1}^6 \alpha_k,
    \\
    \sum_{k=1}^6\frac{\partial}{\partial \alpha_k}
    \left(
    -\frac{1}{2}\sum_{i=1}^4(\tau_i-\pi)^2
    \right)
    &= -\frac{3}{2}\sum_{k=1}^6 \alpha_k + 6\pi,
    \\
    \sum_{k=1}^6\frac{\partial}{\partial \alpha_k}
    \left(
    -\sum_{i=1}^4(\xi-\tau_i)^2
    \right)
    &= -3\sum_{k=1}^6 \alpha_k + 12\xi, 
    \\
    \sum_{k=1}^6\frac{\partial}{\partial \alpha_k}
    \left(
    -\sum_{j=1}^3(\eta_j - \xi)^2
    \right)
    &= -4\sum_{k=1}^6 \alpha_k + 12\xi, 
    \\
    \sum_{k=1}^6\frac{\partial}{\partial \alpha_k}
    \left(
    -\frac{1}{2}\sum_{i=1}^4\sum_{j=1}^3\mathit{Li}_2\left(e^{2\sqrt{-1}(\eta_j-\tau_i)}\right)
    \right)
    &= \frac{\sqrt{-1}}{2}\sum_{i=1}^4\sum_{j=1}^3\log\left(1-e^{2\sqrt{-1}(\eta_j-\tau_i)}\right),
    \\
    \sum_{k=1}^6\frac{\partial}{\partial \alpha_k}
    \left(
    \frac{1}{2}\mathit{Li}_2\left(e^{2\sqrt{-1}(\tau_i-\pi)}\right)
    \right)
    &= -\frac{3\sqrt{-1}}{2}\sum_{i=1}^4 \log\left( 1-e^{2\sqrt{-1}\tau_i} 
    \right),
    \\
    \sum_{k=1}^6\frac{\partial}{\partial \alpha_k}
    \left(
    \sum_{i=1}^4\mathit{Li}_2\left( e^{2\sqrt{-1}(\xi-\tau_i)} \right)
    \right)
    &= 3\sqrt{-1}\sum_{i=1}^4\log\left(1-e^{2\sqrt{-1}(\xi-\tau_i)} \right),
    \\
    \sum_{k=1}^6\frac{\partial}{\partial \alpha_k}
    \left(
    \sum_{j=1}^3 \mathit{Li}_2\left(e^{2\sqrt{-1}(\eta_j -\xi)}\right)
    \right)
    &= -4\sqrt{-1}\sum_{j=1}^3\log\left(1-e^{2\sqrt{-1}(\eta_j-\xi)} \right).
    \\
\end{align*}
This implies 
    \begin{align}\label{kappapf}
    \sum_{k=1}^6 \frac{\partial U}{\partial \alpha_k} 
    =&\  \left[ -2\sum_{k=1}^6\alpha_k + 6\xi + 3\pi  + \frac{\sqrt{-1}}{2}\sum_{i=1}^4\sum_{j=1}^3\log\left(1-e^{2\sqrt{-1}(\eta_j-\tau_i)}\right)
    - \frac{3\sqrt{-1}}{2}\sum_{i=1}^4\log\left(1-e^{2\sqrt{-1}\tau_i}\right) \right. \notag
    \\
    &\ \left. +3\sqrt{-1}\log\left(1-e^{2\sqrt{-1}\xi}\right)
    -\sqrt{-1}\sum_{j=1}^3\log\left(1-e^{2\sqrt{-1}(\eta_j-\xi)}\right)\right]
    -\frac{3}{2} \frac{\partial U}{\partial \xi}.
\end{align}
In particular, by using \eqref{kappa}, \eqref{kappapf} and the fact that $\sum_{i=1}^4 \tau_i = \sum_{k=1}^6 \alpha_k$, when $\partial U/\partial \xi =0$, we have
\begin{align}\label{kappapf2}
\kappa(\alpha_1,\dots,\alpha_6,\xi)=&-\frac{\sqrt{-1}}{2}\sum_{k=1}^6\frac{\partial U}{\partial \alpha_{k}}-\frac{\sqrt{-1}}{2}\sum_{k=1}^6\alpha_{k} + 2\sqrt{-1}\xi-\frac{1}{2}\sum_{i=1}^4\log\big(1-e^{2\sqrt{-1}(\xi-\tau_{i})}\big).
\end{align}

Let $u_i=e^{\sqrt{-1}\alpha_i}$ for $i=1,\dots,6$ and let  $z=e^{-2\sqrt{-1}\xi}.$ 
For a fixed $\alpha$ so that $\mathrm{Re}(\alpha)$ is of the hyperideal type, consider the function $U_\alpha$ of $\xi$ defined by $U_\alpha(\xi)=U(\alpha,\xi).$ From (\ref{partialUxi}), if $\xi$ is a critical point of $U_\alpha,$ then as a necessary condition
\begin{equation}\label{ece}
\frac{(1-z)(1-zu_1u_2u_4u_5)(1-zu_1u_3u_4u_6)(1-zu_2u_3u_5u_6)}{(1-zu_1u_2u_3)(1-zu_1u_5u_6)(1-zu_2u_4u_6)(1-zu_3u_4u_5)}=1,
\end{equation}
which is equivalent to the following quadratic equation 
\begin{equation}\label{qe}
Az^2+Bz+C=0,
\end{equation}
where
\begin{equation*}
\begin{split}
A=&u_1u_4+u_2u_5+u_3u_6-u_1u_2u_6-u_1u_3u_5-u_2u_3u_4-u_4u_5u_6+u_1u_2u_3u_4u_5u_6,\\
B=&-\Big(u_1-\frac{1}{u_1}\Big)\Big(u_4-\frac{1}{u_4}\Big)-\Big(u_2-\frac{1}{u_2}\Big)\Big(u_5-\frac{1}{u_5}\Big)-\Big(u_3-\frac{1}{u_3}\Big)\Big(u_6-\frac{1}{u_6}\Big),\\
C=&\frac{1}{u_1u_4}+\frac{1}{u_2u_5}+\frac{1}{u_3u_6}-\frac{1}{u_1u_2u_6}-\frac{1}{u_1u_3u_5}-\frac{1}{u_2u_3u_4}-\frac{1}{u_4u_5u_6}+\frac{1}{u_1u_2u_3u_4u_5u_6}.
\end{split}
\end{equation*}
Let $\xi(\alpha)$ be the complex number with $\mathrm{Re}(\xi(\alpha))\in [\pi, 2\pi]$ such that 
\begin{equation}\label{za}
e^{-2\sqrt{-1}\xi(\alpha)}=z_{\alpha}=\frac{-B+\sqrt{B^2-4AC}}{2A},
\end{equation}
and let
\begin{equation}\label{za'}
z_{\alpha}'=\frac{-B-\sqrt{B^2-4AC}}{2A}
\end{equation}
be the other root of (\ref{qe}).
Since the region of $\alpha$ is simply connected, we can choose the branch of $\sqrt{B^2-4AC}$ by analytic continuation. A direct computation shows that
\begin{equation}\label{detG}
B^2-4AC=16\det \left[
\begin{array}{cccc}
1 & \cos\alpha_{1} & \cos\alpha_{2} & \cos\alpha_{6}\\
\cos\alpha_{1}& 1 &\cos\alpha_{3} & \cos\alpha_{5}\\
\cos\alpha_{2} & \cos\alpha_{3} & 1 & \cos\alpha_{4} \\
\cos\alpha_{6} & \cos\alpha_{5} & \cos\alpha_{4}  & 1 \\
 \end{array}\right].
\end{equation}
We notice that for a deeply truncated tetrahedron $\Delta(\theta_I,\theta_J)$ with the lengths $\{l_i\}_{i\in I}$ of the edges of deep truncation and the dihedral angles $\{\theta\}_{j\in J}$ at the regular edges, if $\alpha_i=\pi\pm \sqrt{-1}l_i$ for $i\in I$ and $\alpha_j=\pi\pm\theta_j$ for $j\in J,$ then the matrix in (\ref{detG}) equals the Gram matrix of $\Delta(\theta_I,\theta_J)$ as defined in  (\ref{gram2}), which also equals  the value $\mathbb G\big((l_i)_{i\in I},(\theta_j)_{j\in J})$ of the Gram matrix function $\mathbb G$ defined in (\ref{gram}) at $\big((l_i)_{i\in I},(\theta_j)_{j\in J}).$

It is proved in \cite{BY} that 
\begin{equation}\label{=0}
\frac{d U_\alpha}{d \xi}\bigg|_{\xi(\alpha)}=0
\end{equation}
for any $\alpha=(\alpha_1,\dots,\alpha_6)\in\mathbb C^6$ so that $(\mathrm{Re}(\alpha_1),\dots,\mathrm{Re}(\alpha_6))$ is of the hyperideal type, ie, $\xi(\alpha)$ is a critical point of $U_\alpha.$ At this point, we do not know whether $(\alpha,\xi(\alpha))$ always lies in $\mathrm{B_{H,\mathbb C}}.$ It is proved in \cite{BY} that it does when $\alpha_1,\dots,\alpha_6$ are sufficiently close to $\pi.$

For $\alpha\in\mathbb C^6$ so that $(\alpha,\xi(\alpha))\in\mathrm{B_{H,\mathbb C}},$ we define
\begin{equation}\label{W}
W(\alpha)=U(\alpha,\xi(\alpha)).
\end{equation}
Note that by the Chain Rule and (\ref{=0}), we have for $k\in\{1,\dots, 6\},$ 
\begin{equation}\label{3.5}
\begin{split}
\frac{\partial W}{\partial \alpha_{k}}\bigg|_{\alpha}=&\frac{\partial U}{\partial \alpha_k}\bigg|_{(\alpha,\xi(\alpha))}+\frac{\partial U_\alpha}{\partial \xi}\bigg|_{\xi(\alpha)}  \frac{\partial \xi(\alpha)}{\partial \alpha_{k}}\bigg|_{\alpha}=\frac{\partial U}{\partial \alpha_k}\bigg|_{(\alpha,\xi(\alpha))}.
\end{split}
\end{equation}
Finally, we end this section by recalling the following result from \cite{BY} that relates the function $W$ with the co-volume function of deeply truncated tetrahedron.
\begin{theorem}\label{co-vol}\cite[Theorem 3.5]{BY} For a partition $(I,J)$ of $\{1,\dots,6\}$ and a deeply truncated tetrahedron $\Delta$ of type $(I,J),$ we let $\{l_i\}_{i\in I}$ and $\{\theta_i\}_{i\in I}$ respectively be the lengths of and dihedral angles at the edges of deep truncation,  and let $\{\theta_j\}_{j\in J}$ and  $\{l_j\}_{j\in J}$ respectively be the dihedral angles  at and the lengths of the regular edges. Then 
$$W\big((\pi\pm \sqrt{-1} l_i)_{i\in I},(\pi\pm \theta_j)_{j\in J}\big)=2\pi^2+2\sqrt{-1} \mathrm{Cov}\big((l_i)_{i\in I},(\theta_j)_{j\in J}\big)$$
 where $\mathrm{Cov}$ is the co-volume function defined by
 $$\mathrm{Cov}\big((l_i)_{i\in I},(\theta_j)_{j\in J}\big)=\mathrm{Vol}(\Delta)+\frac{1}{2}\sum_{i\in I}\theta_il_i,$$
which for $i\in I$ satisfies
$$\frac{\partial \mathrm{Cov}}{\partial l_i}=\frac{\theta_i}{2}$$
and  for $j\in J$ satisfies
$$\frac{\partial \mathrm{Cov}}{\partial \theta_j}=-\frac{l_j}{2}.$$
\end{theorem}

\section{Asymptotic expansion of the relative Reshetikhin-Turaev invariants}

In \cite{WY2}, the authors studied the exponential growth rate of the relative Reshetikhin-Turaev invariants of the pair $(M,L),$ where $M$ is a closed oriented $3$-manifold and $L$ is a framed link in $M,$ obtained from a fundamental shadow link complement by doing a change-of-pair operation defined in Section \ref{GRT}, and related it to the volume and the Chern-Simons invariant of the hyperbolic cone metric on $M$ with cone angles determined by the sequence of the colorings. The proof of Theorem \ref{main1} is based on the results in \cite{WY2}. We notice that by \cite{CT}, every closed oriented $3$-manifold $M$ can occurs in such pairs $(M,L)$. Throughout this section, we follow the conventions and normalizations of the relative Reshetikhin-Turaev invariants in \cite{WY2}.

\subsection{Fundamental shadow link complements}\label{fundsl}

If we take $c$ truncated tetrahedra $\Delta_1,\dots, \Delta_c$ and glue them together along the triangles of truncation, we obtain a (possibly non-orientable) handlebody of genus $c+1$ with a link in its boundary consisting of the edges of the truncated tetrahedra. By taking the orientable double (the orientable double
covering with the boundary quotient out by the deck involution) of this handlebody, we obtain a link $L_{\text{FSL}}$ inside $M_c=\#^{c+1}(S^2\times S^1).$ We call a link obtained this way a \emph{fundamental shadow link}, and its complement in $M_c$ a \emph{fundamental shadow link complement}.  Alternatively, to construct a fundamental shadow link complement, we can first  take the double of each tetrahedron along the hexagonal faces and then glue the resulting pieces together along homeomorphisms between the $3$-puncture spheres coming from the double of the triangles of truncation.

The fundamental importance of the family of the fundamental shadow link complements  is by \cite{CT} that any compact oriented $3$-manifold with toroidal or empty boundary can be obtained from a suitable fundamental shadow link complement by doing an integral Dehn-filling to some of the boundary components.

\begin{figure}[htbp]
\centering
\includegraphics[scale=0.5]{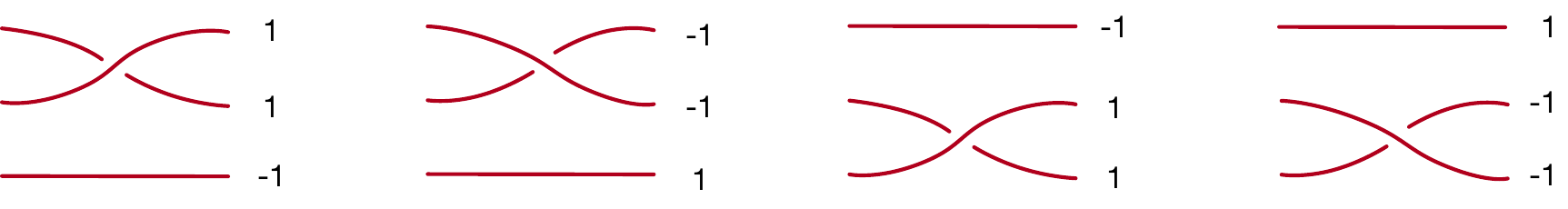}
\caption{ }
\label{mutation}
\end{figure}

A hyperbolic cone metric on $M_c$ with singular locus $L_{\text{FSL}}$ and with sufficiently small cone angles $\theta_1,\dots,\theta_n$ can be constructed as follows. For each $s\in \{1,\dots, c\},$ let $e_{s_1},\dots,e_{s_6}$ be the edges of the building block $\Delta_s,$ and let $\theta_{s_j}$ be half of the cone angle of the component of $L$ containing $e_{s_j}.$ If $\theta_i$'s are sufficiently small, then $\{{\theta_{s_1}},\dots,{\theta_{s_6}}\}$ form the set of dihedral angles of a truncated hyperideal tetrahedron, by abuse of notation still denoted by $\Delta_s.$ Then the hyperbolic cone manifold $M_c$ with singular locus $L_{\text{FSL}}$ and cone angles $\theta_1,\dots,\theta_n$ is obtained by gluing $\Delta_s$'s together along isometries of the triangles of truncation, and taking the double. In this metric, the logarithmic holonomy of the meridian $u_i$ of a tubular neighborhood $N(L_i)$ of $L_i$ satisfies 
\begin{equation}\label{m}
\mathrm{H}(u_i)=\sqrt{-1}\theta_i.
\end{equation}
A preferred longitude $v_i$ on the boundary of $N(L_i)$ can be chosen as follows. Recall that a fundamental shadow link is obtained from the double of a set of truncated tetrahedra (along the hexagonal faces) glued together by orientation preserving homeomorphisms between the trice-punctured spheres coming from the double of the triangles of truncation, and recall also that 
 the mapping class group of trice-punctured sphere is generated by mutations, which could be represented by the four $3$-braids in Figure \ref{mutation}. For each mutation, we assign an integer $\pm 1$ to each component of the braid as in Figure \ref{mutation}; and for a composition of a sequence of mutations, we assign the sum of the  $\pm 1$ assigned by the mutations to each component of the $3$-braid.
In this way, each  orientation preserving homeomorphisms between the trice-punctured spheres assigns three integers to three of the components of $L_{\text{FSL}},$ one for each. For each $i\in \{1,\dots, n\},$ let $\iota_i$ be the sum of all the integers on $L_i$ assigned by the homeomorphisms between the trice-punctured spheres. Then we can choose a preferred longitude $v_i$ on $N(L_i)$ such that $u_i\cdot v_i=1$ and the logarithmic holonomy satisfies
\begin{equation}\label{l}
\mathrm{H}(v_i)=-l_i+\frac{\iota_i\sqrt{-1}\theta_i}{2},
\end{equation}
where $l_i$ is the length of the closed geodesic $L_i.$

\subsection{Growth rate of the relative Reshetikhin-Turaev invariants}\label{GRT}

\begin{figure}[htbp]
\centering
\includegraphics[scale=0.3]{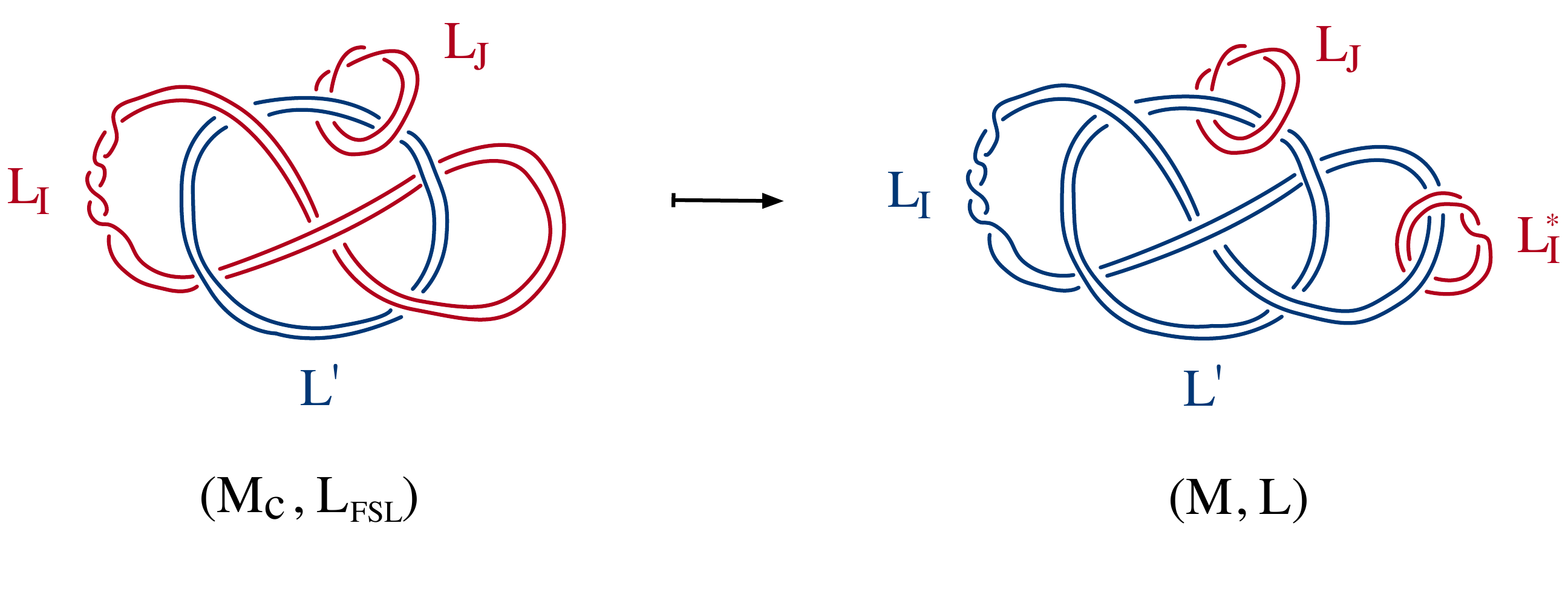}
\caption{Here $L'\subset S^3$ is the disjoint union of $c+1$ unknots with the $0$-framings by doing surgery along which we get $M_c.$ 
}
\label{TFT}
\end{figure}

In this section, we recall the results from \cite{WY2} on the exponential growth rate of the relative Reshetikhin-Turaev invariants of the pair $(M,L)$ obtained from the pair $(M_c,L_{\text{FSL}})$ by doing a change-of-pair operation. More precisely, let $L_1,\dots, L_n$ be the components of $L_{\text{FSL}}.$ For a partition  $(I,J)$ of $\{1,\dots,n\},$ let $L_I=\cup_{i\in I}L_i$ and let $L_J=\cup_{j\in J}L_j.$ Then the pair $(M,L)$ is obtained from $(M_c,L_{\text{FSL}})$ by doing a change-of-pair operation if $M=(M_c)_{L_I}$ is the manifold obtained from $M_c$ by doing a surgery along $L_I$ and $L=L^*_I\cup L_J,$ where $L^*_I=\cup_{i\in I}L^*_i$ and $L^*_i$ is a framed unknot in $M_c\setminus L_{\text{FSL}}$ with the core curve isotopic to the meridian of the tubular neighborhood of $L_i.$ See Figure \ref{TFT}.

By the Symmetric Principle \cite[Lemma 6.3 (iii)]{BHMV}, we can consider the $\mathrm{SO}(3)$-version of the relative Reshetikhin-Turaev invariants, and the $\mathrm{SU}(2)$-version will differ from the $\mathrm{SO}(3)$-version by a nonzero scalar independent of $r.$ 
Let $r$ be an odd integer. Suppose $L^*_i$ has the framing $q_i$ and $L_i$ has the framing  $p_i$  for each $i\in I,$ and  $L_j$ has the framing $p_j$ for each $j\in J.$
Let $(\mathbf b_I^{(r)}, \mathbf a_J^{(r)})$ be a sequence of colorings of $L=L^*_I\cup L_J$, where $b_i^{(r)}, a_j^{(r)} \in \{0,2,\dots, r-3\}^{|L|}$ are even integers for $i\in I, j\in J$. Assume that for any $i\in I,$ either $b_i^{(r)}>\frac{r}{2}$ for all $r$ or $b_i^{(r)}<\frac{r}{2}$ for all $r;$ and for any $j\in J,$ either $a_j^{(r)}>\frac{r}{2}$ for all $r$ or $a_j^{(r)}<\frac{r}{2}$ for all $r.$ In the former case we let $\mu_i=\mu_j=1$ and in the latter case we let $\mu_i=\mu_j=-1,$ and for $i\in I$ we let
$$\theta_i^{(r)}=\mu_i\bigg(\frac{4\pi b_i^{(r)}}{r}-2\pi\bigg)=\bigg|\frac{4\pi b_i^{(r)}}{r}-2\pi\bigg|$$
and for $j\in J$ we let 
$$\theta_j^{(r)}=\mu_j\bigg(\frac{4\pi a_j^{(r)}}{r}-2\pi\bigg)=\bigg|\frac{4\pi a_j^{(r)}}{r}-2\pi\bigg|.$$ 
For simplicity, in the rest of this section we will drop the superscript  and write $\theta_i=\theta_i^{(r)}$ and $\theta_j=\theta_j^{(r)},$ and keep in mind that it is a quantity depending on $r.$  Let $\theta=(\theta_1,\dots,\theta_n),$ and let $M_{L_\theta}$ be the hyperbolic cone manifold consisting of $M$ and a hyperbolic cone metric on $M$  with singular locus $L$ and cone angles $\theta,$ if exists, and let $\mathrm {Vol}(M_{L_\theta})$ and $\mathrm {CS}(M_{L_\theta})$ respectively be the volume and the Chern-Simons invariant of $M_{L_\theta}.$

\begin{theorem}\cite[Theorem 3.1 and Propositions 3.6, 3.7, 5.2, 5.5, 5.9 and 5.10]{WY2}\label{rrt} Suppose $\theta_1,\dots, \theta_n$ converge to sufficiently small limits. Then as $r$ varies over all positive odd integers and at $q=e^{\frac{2\pi\sqrt{-1}}{r}},$ 
\begin{equation}\label{RTV}
\begin{split}
\mathrm{RT}_r&(M,L,(\mathbf b^{(r)}_I,\mathbf a^{(r)}_J))\\
=&\kappa_r \Bigg( \sum_{\epsilon_I\in\{1,-1\}^I}\frac{C^{\epsilon_I}(z^{\epsilon_I})}{\sqrt{-\det\mathrm{Hess}\frac{\mathcal W^{\epsilon_I}(z^{\epsilon_I})}{4\pi\sqrt{-1}}}}\Bigg)e^{\frac{r}{4\pi}\big(\mathrm{Vol}(M_{L_\theta})+\sqrt{-1}\mathrm{CS}(M_{L_\theta})\big)}\Big( 1 + O \Big( \frac{1}{r} \Big) \Big),
\end{split}
\end{equation}
where $\kappa_r,$ $\mathcal W^{\epsilon_I},$ $z^{\epsilon_I}$ and $C^{\epsilon_I}$ are quantities depending on $r$ which will be explained as follows.
\end{theorem}

The quantity $\kappa_r$ is given by
\begin{equation*}
\begin{split}
\kappa_r=&\frac{2^{|I|-2c}}{\{1\}^{|I|-c}}\bigg(\frac{\sin\frac{2\pi}{r}}{\sqrt{r}}\bigg)^{|I|-c}e^{\big(-\sigma(L'\cup L_I)(-\frac{3}{r}-\frac{r+1}{4})-\frac{r}{4}(\sum_{i\in I}q_i+\sum_{k=1}^np_k+2|I|)\big)\sqrt{-1}\pi}\\
=&\frac{(-1)^{-\frac{|I|-c}{2}+\big(-\sigma(L'\cup L_I)(-\frac{3}{r}-\frac{r+1}{4})-\frac{r}{4}(\sum_{i\in I}q_i+\sum_{k=1}^np_k+2|I|)\big)}}{2^{c}}r^{-\frac{|I|-c}{2}},
\end{split}
\end{equation*}
where $\sigma(L'\cup L_I)$ is the signature of the linking matrix of $L'\cup L_I.$

For the function $\mathcal W^{\epsilon_I},$ we let $\alpha_i=\frac{2\pi a_i}{r}$ and $\beta_i=\frac{2\pi b_i^{(r)}}{r}$ for $i\in I,$ $\alpha_j=\frac{2\pi a_j^{(r)}}{r}$ for $j\in J,$ $\xi_s=\frac{2\pi k_s}{r}$ for $s\in\{1,\dots,c\},$ $\tau_{s_i}=\frac{2\pi T_{s_i}}{r}$ for $i\in\{1,\dots,4\},$ and $\eta_{s_j}=\frac{2\pi Q_{s_j}}{r}$ for $j\in\{1,2,3\}.$ For a fixed $(\alpha_j)_{j\in J},$ let
$$\mathrm {D_A}=\Big\{(\alpha_I,\xi)\in\mathbb R^{|I|+c}\ \Big|\ (\alpha_{s_1},\dots,\alpha_{s_6}) \text{ is admissible, } \max\{\tau_{s_i}\}\leqslant \xi_s\leqslant \min\{\eta_{s_j}, 2\pi\}, s\in\{1,\dots,c\}\Big\},$$
where $\alpha_I=(\alpha_i)_{i\in I}$ and $\xi=(\xi_1,\dots,\xi_c),$ and let
$$\mathrm {D_H}=\Big\{(\alpha_I,\xi)\in\mathrm {D_A} \ \Big|\ (\alpha_{s_1},\dots,\alpha_{s_6}) \text{ is of the hyperideal type, } s\in\{1,\dots, c\} \Big\}.$$
For a sufficiently small $\delta >0,$ let 
$$\mathrm {D_H^\delta}=\Big\{(\alpha_I,\xi)\in\mathrm {D_H}\ \Big|\ d((\alpha_I,\xi), \partial\mathrm {D_H})>\delta \Big\},$$
where $d$ is the Euclidean distance on $\mathbb R^n.$ Then $\mathcal W^{\epsilon_I}$ is the following function
\begin{equation}\label{WforRT}
\begin{split}
\mathcal W^{\epsilon_I}(\alpha_I,\xi)=&-\sum_{i\in I}q_i(\beta_i-\pi)^2-\sum_{j\in J}p_j(\alpha_j-\pi)^2\\
&-\sum_{i\in I}p_i(\alpha_i-\pi)^2-\sum_{i\in I}2\epsilon_i(\alpha_i-\pi)(\beta_i-\pi)\\
&-\sum_{i=1}^n\frac{\iota_i}{2}(\alpha_i-\pi)^2+\sum_{s=1}^c U(\alpha_{s_1},\dots,\alpha_{s_6},\xi_s)+\Big(\sum_{i=1}^n\frac{\iota_i}{2}\Big)\pi^2
\end{split}
\end{equation}
where $U$ is as defined in (\ref{term}), which is continuous on 
$$\mathrm{D_{H,\mathbb C}}=\big\{(\alpha_I,\xi_1)\in\mathbb C^{|I|+c}\ \big|\ (\mathrm{Re}(\alpha_I),\mathrm{Re}(\xi))\in \mathrm{D_{H}}\big\}$$ and for any $\delta>0$ is analytic on 
$$\mathrm{D^\delta_{H,\mathbb C}}=\big\{(\alpha_I,\xi_1)\in\mathbb C^{|I|+c}\ \big|\ (\mathrm{Re}(\alpha_I),\mathrm{Re}(\xi))\in \mathrm{D^\delta_{H}}\big\},$$ 
where $\mathrm{Re}(\alpha_I)=(\mathrm{Re}(\alpha_i))_{i\in I}$ and $\mathrm{Re}(\xi)=(\mathrm{Re}(\xi_1),\dots, \mathrm{Re}(\xi_c)).$

For $z^{\epsilon_I},$ we for each $i\in I,$ let $\mathrm H(u_i)$ be the logarithmic holonomy of $u_i$ in the hyperbolic cone manifold $M_{L_\theta}$ and let
$$\alpha^*_i=\pi+\frac{\epsilon_i\mu_i\sqrt{-1}}{2}\mathrm H(u_i).$$
For $s\in\{1,\dots,s\},$ let $I_s=\{s_1,\dots,s_6\}\cap I$ and let $J_s=\{s_1,\dots,s_6\}\cap J.$ Let $\alpha_{I_s}^*=(\alpha^*_{s_i})_{s_i\in I_s},$ $\alpha_{J_s}=(\alpha_{s_j})_{j\in J}$ and 
$$\xi_s^*=\xi(\alpha_{I_s}^*,\alpha_{J_s}).$$
Suppose $\theta_1,\dots,\theta_n$ are sufficiently small. Then by \cite[Proposition 5.2]{WY2}, $\mathcal W^{\epsilon_I}$ has a critical point 
$$z^{\epsilon_I}=\Big(\big(\alpha^*_i\big)_{i\in I}, \big(\xi^*_s\big)_{s=1}^c\Big)$$
in $\mathrm{D_{H,\mathbb C}}$  with critical value 
$$2c\pi^2+\sqrt{-1}\Big(\mathrm{Vol}(M_{L_\theta})+\sqrt{-1}\mathrm{CS}(M_{L_\theta})\Big).$$

By the computation in the end of the proof of \cite[Proposition 5.5]{WY2}, we have
\begin{equation}\label{CRT}
\begin{split}
C^{\epsilon_I}=\frac{(-1)^{-\frac{rc}{2}}r^{\frac{|I|-c}{2}}}{2^{\frac{3|I|-c}{2}}\pi^{\frac{|I|+c}{2}}}e^{\sqrt{-1}\big(\sum_{i\in I}q_i\beta_i+\sum_{i\in I}(p_i+\frac{\iota_i}{2})\alpha^*_i+\sum_{j\in J}(p_j+\frac{\iota_j}{2})\alpha_j+\sum_{i\in I}\epsilon_i(\alpha^*_i+\beta_i)\big)+\sum_{s=1}^c\kappa(\alpha_{I_s}^*,\alpha_{J_s},\xi_s^*)},
\end{split}
\end{equation}
where $\kappa$ is as defined in (\ref{kappa}).

Finally, we notice that all of $\mathcal W^\epsilon,$ $z^\epsilon$ and $C^\epsilon$ depend on $r$ because each $\beta_i=\frac{2\pi b_i^{(r)}}{r},$ $i\in I,$ and $\alpha_j=\frac{2\pi a_j^{(r)}}{r},$ $j\in J,$  does.


\subsection{Proof of Theorem \ref{main1}}

Following the notations in Sections \ref{fundsl} and \ref{GRT}, for $k\in \{1,\dots, n\},$ let $u_k$ and $v_k$ respectively be the meridian and the preferred longitude of a tubular neighborhood of $L_k.$ For $i\in I,$ let $m_i=p_iu_i+v_i$ be the meridian of a tubular neighborhood of $L_i^*$ and let $\mathbf m=\big((m_i)_{i\in I},(u_j)_{j\in J}\big).$ For each $i\in I,$ let $\gamma_i=-u_i+q_i(p_iu_i+v_i)=-u_i+q_im_i$ be the parallel copy of $L^*_i$ that is isotopic to $L^*_i$ given by the framing $q_i$ of $L^*_i$ and with the orientation so that $m_i \cdot \gamma_i=1;$ and for each $j\in J,$ let $\gamma_j=p_ju_j+v_j$ be the parallel copy of $L_j$ that is isotopic to $L_j$ given by the framing $p_j$ of $L_j$ and with the orientation so that $u_j\cdot \gamma_j=1.$ Let $\rho_{M_{L_\theta}}:\pi_1(M\setminus L)\to\mathrm{PSL}(2;\mathbb C)$ be the holonomy representation of the restriction of hyperbolic cone metric $M_{L_\theta}$ to $M\setminus L,$ and let  $\mathbb{T}_{(M\setminus L,\mathbf m)}([\rho_{M_{L_\theta}}])$ be the twisted Reidemeister torsion of $M\setminus L$ twisted by the adjoint action of $\rho_{M_{L_\theta}}.$ For each $k\in \{1,\dots,n\},$ let $\mathrm H(\gamma_k)$ be the logarithmic holonomy of $\gamma_k$ in $\rho_{M_{L_\theta}}.$ Then Theorem \ref{main1} can be rephrased as the following Theorem \ref{main1'}.

\begin{theorem}\label{main1'} Suppose $\theta_1,\dots, \theta_n$ converge to sufficiently small limits. Then as $r$ varies over all positive odd integers and at $q=e^{\frac{2\pi \sqrt{-1}}{r}},$ 
$$\mathrm{RT}_r(M,L,(\mathbf b^{(r)}_I,\mathbf a^{(r)}_J))=C \frac{e^{\frac{1}{2}\sum_{k=1}^n\mu_k\mathrm H(\gamma_k)}}{\sqrt{\pm\mathbb{T}_{(M\setminus L,\mathbf m)}([\rho_{M_{L_\theta}}])}}e^{\frac{r}{4\pi}\big(\mathrm{Vol}(M_{L_\theta})+\sqrt{-1}\mathrm{CS}(M_{L_\theta})\big)}\Big( 1 + O \Big( \frac{1}{r} \Big) \Big),$$
where 
$$C={(-1)^{\sum_{i\in I} \big(p_i+\frac{\iota_i}{2}+q_i\big)+\sum_{j\in J}\big(p_j+\frac{\iota_j}{2}\big)-\sigma(L'\cup L_I)\big(-\frac{3}{r}-\frac{r+1}{4}\big)-\frac{r}{4}\big(\sum_{i\in I}q_i+\sum_{k=1}^np_k\big)+\frac{|I|}{2}+(|I|-c)\big(\frac{r}{2}-\frac{1}{4}\big)}}$$
is a quantity of norm $1,$ which is independent of the geometric structure on $M.$ 
\end{theorem}

By Theorem \ref{rrt}, to prove Theorem \ref{main1'}, we need to compute the $ \frac{C^{\epsilon_I}}{\sqrt{-\det\mathrm{Hess}\frac{\mathcal W^{\epsilon_I}(z^{\epsilon_I})}{4\pi\sqrt{-1}}}}$ for each $\epsilon_I\in\{1,-1\}^I.$ To do this, we need the following Lemmas \ref{L3.1}, \ref{L3.2}, \ref{L3.4} and \ref{L3.5}.

For  each $s\in\{1,\dots,c\},$ we let $I_s=\{s_1,\dots,s_6\}\cap I$ and $J_s=\{s_1,\dots,s_6\}\cap J.$ We also let $\alpha_{I_s}^*=(\alpha^*_{s_i})_{s_i\in I_s},$ $\alpha_{J_s}=(\alpha_{s_j})_{j\in J},$ $\xi_s^*=\xi(\alpha_{I_s}^*,\alpha_{J_s})$ and $z_s^*=(\alpha_{I_s}^*,\alpha_{J_s},\xi_s^*).$ 

\begin{lemma}\label{L3.1}  For each $i\in I,$ recall that $m_i=p_iu_i+v_i$ is the meridian of a tubular neighborhood of $L^*_i.$ Then 
$$-\det\mathrm{Hess}\mathcal W^{\epsilon_I}(z^{\epsilon_I})=-(-2)^{|I|}\det\bigg(\frac{\partial \mathrm H(m_{i_1})}{\partial \mathrm H(u_{i_2})}\bigg)_{i_1,i_2\in I} \prod _{s=1}^{c} \frac{\partial ^2U}{\partial \xi_s^2}\bigg|_{z^*_s}.$$
\end{lemma}

\begin{proof} For $s\in\{1,\dots,c\}$ and $i\in I,$ we denote by $s\sim i$ if the tetrahedron $\Delta_s$ intersects the component $L_i$ of $L_{\text{FSL}},$  and for $\{i_1, i_2\}\subset I$ we denote by $s\sim i_1,i_2$ if $\Delta_s$ intersects both $L_{i_1}$ and $L_{i_2}.$ For $s\in\{1,\dots,c\},$ let $\alpha_s=(\alpha_{s_1},\dots,\alpha_{s_6})$ and let $\alpha_s^*=(\alpha_{I_s}^*,\alpha_{J_s}).$ Then we have the following claims:
\begin{enumerate}[(1)]
\item For $s\in\{1,\dots,c\},$
$$\frac{\partial^2 \mathcal W^{\epsilon_I}}{\partial \xi_s^2}\bigg|_{z^{ \epsilon_I}}=\frac{\partial^2 U}{\partial \xi_s^2}\bigg|_{z_s^*}.$$

\item For $\{s_1,s_2\}\subset\{1,\dots,c\},$
$$\frac{\partial^2 \mathcal W^{ \epsilon_I}}{\partial \xi_{s_1}\partial \xi_{s_2}}\bigg|_{z^{ \epsilon_I}}=0.$$

\item For $i\in I$ and $s\in\{1,\dots,c\},$
$$\frac{\partial^2 \mathcal W^{ \epsilon_I}}{\partial \alpha_i\partial \xi_s}\bigg|_{z^{ \epsilon_I}}=-\frac{\partial^2 U}{\partial \xi_s^2}\bigg|_{z^*_s} \frac{\partial \xi_s(\alpha_s)}{\partial \alpha_i}\bigg|_{\alpha_s^*}.$$

\item For $i\in I,$ 
$$\frac{\partial^2 \mathcal W^{ \epsilon_I}}{\partial \alpha_i^2}\bigg|_{z^{ \epsilon_I}}=-2\frac{\partial \mathrm H(m_i)}{\partial \mathrm H(u_i)}+\sum_{s\sim i}\frac{\partial^2 U}{\partial \xi_s^2}\bigg|_{z^*_s} \bigg(\frac{\partial \xi_s(\alpha_s)}{\partial \alpha_i}\bigg|_{\alpha_s^*}\bigg)^2.$$

\item For $\{i_1,i_2\}\subset I,$ 
$$\frac{\partial^2 \mathcal W^{ \epsilon_I}}{\partial \alpha_{i_1}\partial \alpha_{i_2}}\bigg|_{z^{ \epsilon_I}}=-2\frac{\epsilon_{i_1}\mu_{i_1}}{\epsilon_{i_2}\mu_{i_2}} \frac{\partial \mathrm H(m_{i_1})}{\partial \mathrm H(u_{i_2})}+\sum_{s\sim i_1,i_2}\frac{\partial^2 U}{\partial \xi_s^2}\bigg|_{z^*_s} \frac{\partial \xi_s(\alpha_s)}{\partial \alpha_{i_1}}\bigg|_{\alpha_s^*} \frac{\partial \xi_s(\alpha_s)}{\partial \alpha_{i_2}}\bigg|_{\alpha_s^*}.$$
\end{enumerate}
We will prove the claims (1) -- (5) in the end. Assuming these claims, then 
\begin{equation}\label{congurent0}
\mathrm{Hess}\mathcal W^{{ \epsilon_I}}(z^{{ \epsilon_I}})=A\cdot D\cdot A^T,
\end{equation}
with $D$ and $A$ defined as follows. 
The matrix $D$ is a block matrix with the left-top block the $|I|\times|I|$ matrix 
$$\bigg(-2\frac{\epsilon_{i_1}\mu_{i_1}}{\epsilon_{i_2}\mu_{i_2}} \frac{\partial \mathrm H(m_{i_1})}{\partial \mathrm H(u_{i_2})}\bigg)_{i_1,i_2\in I},$$
the right-top and the left-bottom blocks consisting of $0$'s, and the right-bottom block the $c \times c$ diagonal matrix with the diagonal entries
$\frac{\partial^2 U}{\partial \xi_1^2}\Big|_{z_1^*},\dots, \frac{\partial^2 U}{\partial \xi_c^2}\Big|_{z_c^*}.$
Then 
\begin{equation}\label{detD0}
\begin{split}
\det D=&\ (-2)^{|I|}\det\bigg(\frac{\epsilon_{i_1}\mu_{i_1}}{\epsilon_{i_2}\mu_{i_2}} \frac{\partial \mathrm H(m_{i_1})}{\partial \mathrm H(u_{i_2})}\bigg)_{i_1,i_2\in I} \prod _{s=1}^{c} \frac{\partial ^2U}{\partial \xi_s^2}\bigg|_{z^*_s}\\
=&\ (-2)^{|I|}\det\bigg(\frac{\partial \mathrm H(m_{i_1})}{\partial \mathrm H(u_{i_2})}\bigg)_{i_1,i_2\in I}  \prod _{s=1}^c \frac{\partial ^2U}{\partial \xi_s^2}\bigg|_{z^*_s}.
\end{split}
\end{equation}
The matrix $A$ is a block matrix with the left-top and the right-bottom blocks respectively the $|I|\times|I|$ and $c\times c$ identity matrices, the left-bottom block consisting of $0$'s and the right-top block the $|I|\times c$ matrix with entries $a_{is},$ $i\in I$ and $s\in\{1,\dots,c\},$ given by
$$a_{is}=-\frac{\partial \xi_s(\alpha_s)}{\partial \alpha_i}\bigg|_{\alpha_s^*}$$
if $s\sim i$ and $a_{is}=0$
if otherwise. Since $A$ is upper triangular with all diagonal entries equal to $1,$ 
\begin{equation}\label{detA0}
\det A=1.
\end{equation}
The result then follows from (\ref{congurent0}), (\ref{detD0}) and (\ref{detA0}), and we are left to prove the claims (1) -- (5). 

Claims (1) and (2) are straightforward  from the definition of $\mathcal W^{ \epsilon_I}.$ For (3), we have
\begin{equation}\label{WU0}
\frac{\partial \mathcal W^{ \epsilon_I}}{\partial \xi_s}\bigg|_{\big((\alpha_i)_{i\in I},\xi_1,\dots,\xi_c\big)}=\frac{\partial U}{\partial \xi_s}\bigg|_{(\alpha_s,\xi_s)}.
\end{equation}
Let 
$$f(\alpha_s,\xi_s)\doteq \frac{\partial U}{\partial \xi_s}\bigg|_{(\alpha_s,\xi_s)}$$
and 
$$g(\alpha_s)\doteq f(\alpha_s,\xi_s(\alpha_s)).$$
Then 
$$g(\alpha_s)=\frac{\partial U}{\partial \xi_s}\bigg|_{(\alpha_s,\xi_s(\alpha_s))}=\frac{dU_{\alpha_s}}{d \xi_s}\bigg|_{\xi_s(\alpha_s)}\equiv 0,$$
and hence
\begin{equation}\label{g=00}
\frac{\partial g}{\partial \alpha_{s_i}}\bigg|_{\alpha_s}=0.
\end{equation}
On the other hand, we have
\begin{equation}\label{+0}
\begin{split}
\frac{\partial g}{\partial \alpha_{s_i}}\bigg|_{\alpha_s}=&\frac{\partial f}{\partial \alpha_{s_i}}\bigg|_{(\alpha_s,\xi_s(\alpha_s))}+\frac{\partial f}{\partial \xi_s}\bigg|_{(\alpha_s,\xi_s(\alpha_s))}  \frac{\partial \xi_s(\alpha_s)}{\partial \alpha_{s_i}}\bigg|_{\alpha_s}\\
=&\frac{\partial^2 U}{\partial \alpha_{s_i}\partial \xi_s}\bigg|_{(\alpha_s,\xi_s(\alpha_s))}+\frac{\partial^2 U}{\partial \xi_s^2}\bigg|_{(\alpha_s,\xi_s(\alpha_s))} \frac{\partial \xi_s(\alpha_s)}{\partial \alpha_{s_i}}\bigg|_{\alpha_s}.
\end{split}
\end{equation}
Putting (\ref{g=00}) and (\ref{+0}) together, we have
\begin{equation}\label{alphaxi0}
\frac{\partial^2 U}{\partial \alpha_{s_i}\partial \xi_s}\bigg|_{(\alpha_s,\xi_s(\alpha_s))}=-\frac{\partial^2 U}{\partial \xi_s^2}\bigg|_{(\alpha_s,\xi_s(\alpha_s))} \frac{\partial \xi_s(\alpha_s)}{\partial \alpha_{s_i}}\bigg|_{\alpha_s},
\end{equation}
and (3) follows from (\ref{WU0}) and (\ref{alphaxi0}).

For (4) and (5), we have 
\begin{equation}\label{3.10}
\frac{\partial ^2\mathcal W^{ \epsilon_I}}{\partial \alpha_i^2}\bigg|_{z^{ \epsilon_I}}=-2p_i-\iota_i+\sum_{s\sim i}\frac{\partial ^2U}{\partial \alpha_i^2}\bigg|_{z^*_s},
\end{equation}
and
\begin{equation}\label{3.11}
\frac{\partial ^2\mathcal W^{ \epsilon_I}}{\partial \alpha_{i_1}\partial \alpha_{i_2}}\bigg|_{z^{ \epsilon_I}}=\sum_{s\sim i_1,i_2}\frac{\partial ^2U}{\partial \alpha_{i_1}\partial \alpha_{i_2}}\bigg|_{z^*_s}.
\end{equation}
By using (\ref{3.5}) and the Chain Rule, for  $j,k\in\{1,\dots,6\}$ we have 
\begin{equation*}\label{3.12}
\begin{split}
\frac{\partial^2 W}{\partial \alpha_{s_j}\partial \alpha_{s_k}}\bigg|_{\alpha_s}
=\frac{\partial ^2U}{\partial \alpha_{s_j}\partial \alpha_{s_k}}\bigg|_{(\alpha_s,\xi_s(\alpha_s))}+\frac{\partial^2 U}{\partial \alpha_{s_k}\partial \xi_s}\bigg|_{(\alpha_s,\xi_s(\alpha_s))}  \frac{\partial \xi_s(\alpha_s)}{\partial \alpha_{s_j}}\bigg|_{\alpha_s}.
\end{split}
\end{equation*}
Together with (\ref{alphaxi0}), for $j,k\in\{1,\dots, 6\}$ we have
\begin{equation}\label{3.13}
\begin{split}
\frac{\partial^2 U}{\partial \alpha_{s_j}\partial \alpha_{s_k}}\bigg|_{(\alpha_s,\xi_s(\alpha_s))}=&\frac{\partial^2 W}{\partial \alpha_{s_j}\partial \alpha_{s_k}}\bigg|_{\alpha_s}-\frac{\partial^2 U}{\partial \alpha_{s_k}\partial \xi_s}\bigg|_{(\alpha_s,\xi_s(\alpha_s))}  \frac{\partial \xi_s(\alpha_s)}{\partial \alpha_{s_j}}\bigg|_{\alpha_s}\\
=&\frac{\partial^2 W}{\partial \alpha_{s_j}\partial \alpha_{s_k}}\bigg|_{\alpha_s}+\frac{\partial^2 U}{\partial \xi_s^2}\bigg|_{(\alpha_s,\xi_s(\alpha_s))} \frac{\partial \xi_s(\alpha_s)}{\partial \alpha_{s_j}}\bigg|_{\alpha_s}  \frac{\partial \xi_s(\alpha_s)}{\partial \alpha_{s_k}}\bigg|_{\alpha_s}.
\end{split}
\end{equation}
By (\ref{3.10}), (\ref{3.11}) and (\ref{3.13}) and we have
\begin{equation}\label{3.15}
\begin{split}
\frac{\partial ^2\mathcal W^{ \epsilon_I}}{\partial \alpha_i^2}\bigg|_{z^{ \epsilon_I}}
=&-2p_i-\iota_i+\sum_{s\sim i}\sum_{s_k}\frac{\partial^2 W}{\partial \alpha_{s_k}^2}\bigg|_{\alpha^*_s}+\sum_{s\sim i}\frac{\partial^2 U}{\partial \xi_s^2}\bigg|_{z^*_s} \bigg(\frac{\partial \xi_s(\alpha_s)}{\partial \alpha_i}\bigg|_{\alpha_s^*}\bigg)^2,
\end{split}
\end{equation}
where the second sum in the third term of the right hand side is over $s_k$ such that the edge $e_{s_k}$ in $\Delta_s$ intersects the component  $L_i;$ and
\begin{equation}\label{3.17}
\begin{split}
\frac{\partial ^2\mathcal W^{ \epsilon_I}}{\partial \alpha_{i_1}\partial \alpha_{i_2}}\bigg|_{z^{ \epsilon_I}}=&\sum_{s\sim i_1,i_2}\sum_{s_j,s_k}\frac{\partial^2 W}{\partial \alpha_{s_j}\partial \alpha_{s_k}}\bigg|_{\alpha^*_s}+\sum_{s\sim i_1,i_2}\frac{\partial^2 U}{\partial \xi_s^2}\bigg|_{z^*_s} \frac{\partial \xi_s(\alpha_s)}{\partial \alpha_{i_1}}\bigg|_{\alpha_s^*} \frac{\partial \xi_s(\alpha_s)}{\partial \alpha_{i_2}}\bigg|_{\alpha_s^*},
\end{split}
\end{equation}
where the second sum in the first term of the right hand side is over $s_j,s_k$ such that the edge $e_{s_j}$ in $\Delta_s$ intersects the component  $L_{i_1}$ and the edge $e_{s_k}$ in $\Delta_s$ intersects the component  $L_{i_2}.$

At a hyperbolic cone metric on $M_c$ with singular locus $L_{\text{FSL}},$ by Theorem \ref{co-vol},  for $i,j\in\{1,\dots, 6\}$ we have
\begin{equation}\label{3.18}
\frac{\partial^2 W}{\partial \alpha_{s_i}\partial \alpha_{s_j}}\bigg|_{\alpha^*_s}=-\sqrt{-1}\frac{\epsilon_{s_i}\mu_{s_i}}{\epsilon_{s_j}\mu_{s_j}} \frac{\partial l_{s_i}}{\partial \theta_{s_j}},
\end{equation}
where $l_{s_k}$ is the length of $e_{s_k}$ of $\Delta_s,$ and  if  $e_{s_k}$ intersects $L_i$ then $\epsilon_{s_k}=\epsilon_i,$ $\mu_{s_k}=\mu_i$ and $\theta_{s_k}=\frac{\theta_i}{2}$ is the half of the cone angle at $L_i.$ We also observe that that
\begin{equation}\label{3.14}
l_i=\sum_{s\sim i}\sum_{s_k}l_{s_k},
\end{equation}
where the second sum is over $s_k$ such that the edge $e_{s_k}$ in $\Delta_s$ intersects the component  $L_i.$  

Then by (\ref{3.18}), (\ref{3.14}) (\ref{m}) and (\ref{l}) we have
\begin{equation}\label{3.16}
\begin{split}
-2p_i-\iota_i+\sum_{s\sim i}\sum_{s_k}\frac{\partial^2 W}{\partial \alpha_{s_k}^2}\bigg|_{\alpha^*_s}=&-2p_i-\iota_i-\sqrt{-1}\sum_{s\sim i}\sum_{s_k}\frac{\partial l_{s_k}}{\partial \theta_{s_k}}\\=&-2p_i-\iota_i-2\sqrt{-1}\frac{\partial l_i}{\partial \theta_i}\\
=&-2\frac{\partial \big(p_i\sqrt{-1}\theta_i-l_i+\frac{\iota_i}{2}\sqrt{-1}\theta_i \big)}{\partial (\sqrt{-1}\theta_i)}\\
=&-2\frac{\partial \big(p_i\mathrm H(u_i) +\mathrm H(v_i)\big)}{\partial (\sqrt{-1}\theta_i)}=-2\frac{\partial \mathrm H(m_i)}{\partial \mathrm H(u_i)}.
\end{split}
\end{equation}
From (\ref{3.15}) and (\ref{3.16}), (4) holds at  hyperbolic cone metrics on $M_c$ with singular locus $L_{\text{FSL}}.$ By the analyticity of the involved functions (see for e.g. \cite[Lemma 4.2]{WY2}),  (\ref{3.16}) still holds in a neighborhood of the complete hyperbolic structure on $M_c\setminus L_{\text{FSL}},$ from which (4) follows. 

By (\ref{3.18}), (\ref{3.14}), (\ref{m}) and (\ref{l}) we have
\begin{equation}\label{3.19}
\begin{split}
\sum_{s\sim i_1,i_2}\sum_{s_j,s_k}\frac{\partial^2 W}{\partial \alpha_{s_j}\partial \alpha_{s_k}}\bigg|_{\alpha^*_s}=&-2\sqrt{-1}\sum_{s\sim i_1}\sum_{s_j}\frac{\epsilon_{s_{j}}\mu_{s_{j}}}{\epsilon_{{i_2}}\mu_{{i_2}}}\frac{\partial l_{s_j}}{\partial \theta_{i_2}}\\
=&-2\sqrt{-1}\frac{\epsilon_{i_1}\mu_{i_1}}{\epsilon_{{i_2}}\mu_{{i_2}}}\frac{\partial l_{i_1}}{\partial \theta_{i_2}}\\
=&-2\frac{\epsilon_{i_1}\mu_{i_1}}{\epsilon_{{i_2}}\mu_{{i_2}}}\frac{\partial \big( p_{i_1}\sqrt{-1}\theta_{i_1}-l_{i_1}+\frac{\iota_{i_1}}{2}\sqrt{-1}\theta_{i_1}\big)}{\partial (\sqrt{-1}\theta_{i_2})}\\
=&-2\frac{\epsilon_{i_1}\mu_{i_1}}{\epsilon_{{i_2}}\mu_{{i_2}}}\frac{\partial \big(p_{i_1}\mathrm H(u_{i_1}) +\mathrm H(v_{i_1})\big)}{\partial (\sqrt{-1}\theta_{i_2})}\\=&-2\frac{\epsilon_{i_1}\mu_{i_1}}{\epsilon_{{i_2}}\mu_{{i_2}}}\frac{\partial \mathrm H(m_{i_1})}{\partial \mathrm H(u_{i_2})},
\end{split}
\end{equation}
where the second sum on the left hand side is over $s_j,s_k$ such that the edge $e_{s_j}$ in $\Delta_s$ intersects the component  $L_{i_1}$ and the edge $e_{s_k}$ in $\Delta_s$ intersects the component $L_{i_2},$ the second sum on the right hand side of the first equation is over $s_j$ such that the edge $e_{s_j}$ in $\Delta_s$ intersects the component  $L_{i_1},$ and the third equality comes from the fact that
$$\frac{\partial (\sqrt{-1}\theta_{i_1})}{\partial (\sqrt{-1}\theta_{i_2})}=\frac{\partial \mathrm H(u_{i_1})}{\partial \mathrm H(u_{i_2})}=0.$$
 From (\ref{3.17}) and (\ref{3.19}), (5) holds at  hyperbolic cone metrics on $M_c\setminus L_{\text{FSL}}.$ By the analyticity of the involved functions, (\ref{3.19}) still holds in a neighborhood of the complete hyperbolic structure on $M_c\setminus L_{\text{FSL}},$ from which (5) follows. 
\end{proof}

\begin{lemma}\label{L3.2} 
\begin{equation*}
\begin{split}
\kappa(\alpha_{I_s}^*,\alpha_{J_s},\xi_s^*)=&-\frac{\sqrt{-1}}{2}\sum_{k=1}^6\frac{\partial U}{\partial \alpha_{s_k}}\bigg|_{z^*_s}\\
&-\frac{\sqrt{-1}}{2}\sum_{s_i\in I_s}\alpha^*_{s_i}-\frac{\sqrt{-1}}{2}\sum_{s_j\in J_s}\alpha_{s_j}+2\sqrt{-1}\xi_s^*-\frac{1}{2}\sum_{i=1}^4\log\big(1-e^{2\sqrt{-1}(\xi_s^*-\tau^*_{s_i})}\big),
\end{split}
\end{equation*}
where with the notation that $\alpha^*_{s_j}=\alpha_{s_j}$ for $s_j\in J_s,$
$\tau^*_{s_1}=\frac{\alpha^*_{s_1}+\alpha^*_{s_2}+\alpha^*_{s_3}}{2},$ $\tau^*_{s_2}=\frac{\alpha^*_{s_1}+\alpha^*_{s_5}+\alpha^*_{s_6}}{2},$ $\tau^*_{s_3}=\frac{\alpha^*_{s_2}+\alpha^*_{s_4}+\alpha^*_{s_6}}{2}$ and  $\tau^*_{s_4}=\frac{\alpha^*_{s_3}+\alpha^*_{s_4}+\alpha^*_{s_5}}{2}.$ 
\end{lemma}

\begin{proof}
  The result follows directly from \eqref{kappapf2}.
\end{proof}

\begin{lemma}\label{L3.4} For $i\in I$ recall that $\gamma_i=-u_i+q_i(p_iu_i+v_i)$ is the parallel of copy of $L^*_i$ given by the framing $q_i,$ and for each $j\in J$ recall that $\gamma_j=p_ju_j+v_j$ is the parallel copy of $L_j$ given by the framing $p_j.$  Then
\begin{equation*}
\begin{split}
&\sqrt{-1}\bigg(\sum_{i\in I}q_i\beta_i+\sum_{i\in I}(p_i+\frac{\iota_i}{2})\alpha^*_i+\sum_{j\in J}(p_j+\frac{\iota_j}{2})\alpha_j+\sum_{i\in I}\epsilon_i(\alpha^*_i+\beta_i)\bigg)-\frac{\sqrt{-1}}{2}\sum_{s=1}^c\bigg(\sum_{k=1}^6\frac{\partial U}{\partial \alpha_{s_k}}\bigg|_{z_s^*}\bigg)\\
=&\bigg(\sum_{i\in I} \Big(p_i+\frac{\iota_i}{2}+q_i+2\epsilon_i\Big)+\sum_{j\in J}\Big(p_j+\frac{\iota_j}{2}\Big)\bigg)\sqrt{-1}\pi+\frac{1}{2}\sum_{k=1}^n\mu_k\mathrm H(\gamma_k).
\end{split}
\end{equation*}
\end{lemma}

\begin{proof} Note that both sides of the equality depend analytically on $(\mathrm{H}(u_i))_{i\in I}$. By \cite[Lemma 4.2]{WY2}, it suffices to prove the result for the special case that $\mathrm{H}(u_i) = \sqrt{-1}\theta_i$ for $i\in I$. Note that geometrically this corresponds to the case where $(M,L)=(M_c,L_{\text{FSL}})$ with the hyperbolic cone manifold $M_c\setminus L_{\text{FSL}}$ obtained by gluing hyperideal tetrahedra $\Delta_1,\dots,\Delta_s$ together along the truncated faces and then taking the orientable double. In particular, we have
\begin{align}\label{Sect3HuHv}
    \mathrm H(u_i)=\sqrt{-1}\theta_i , \quad 
\mathrm H(v_i)=-l_i+\frac{\iota_i}{2}\sqrt{-1}\theta_i.
\end{align}

To prove the result in this special case, we first show that
\begin{equation}\label{L3.3}
-\frac{\sqrt{-1}}{2}\sum_{s=1}^c\bigg(\sum_{k=1}^6\frac{\partial U}{\partial \alpha_{s_k}}\bigg|_{z_s^*}\bigg)=-\frac{1}{2}\sum_{i\in I}\epsilon_i\mu_i\Big(\mathrm H(v_i)-\frac{\iota_i}{2}\mathrm H(u_i)\Big)-\frac{1}{2}\sum_{j\in J}\mu_jl_j.
\end{equation}
For each $s\in \{1,\dots, c\}$ let $e_{s_1},\dots,e_{s_6}$ be the edges of $\Delta_s$ and for each $k\in\{1,\dots,6\}$ let $l_{s_k}$ and $\theta_{s_k}$ respectively be the length of and the dihedral angle at $e_{s_k}.$ If $e_{s_k}$ intersects the component $L_i$ of $L_{\text{FSL}}$ for some $i\in I,$ then $\mathrm H(u_i)=\sqrt{-1} \theta_i=2\sqrt{-1}\theta_{s_k}$ and let $\alpha_{s_k}=\alpha^*_i=\pi+\frac{\epsilon_i\mu_i\sqrt{-1}\mathrm H(u_i)}{2}=\pi-\epsilon_{s_k}\mu_{s_k}\theta_{s_k},$ where $\epsilon_{s_k}=\epsilon_i$ and $\mu_{s_k}=\mu_i;$ and if $e_{s_k}$ intersects the component $L_j$ of $L_{\text{FSL}}$ for some $j\in J,$ then $\theta_j=2\theta_{s_k}$ and let $\alpha_{s_k}=\alpha_j=\pi+\frac{\mu_j\theta_j}{2}=\pi+\mu_{s_k}\theta_{s_k},$ where $\mu_{s_k}=\mu_j.$ By \eqref{3.5} and Theorem \ref{co-vol}, we have for $s_k\in I_s$
\begin{equation}\label{3.3}
\frac{\partial U}{\partial \alpha_{s_k}}\bigg|_{z_s^*}
=
\frac{\partial W}{\partial \alpha_{s_k}}\bigg|_{\big(a_{I_s}^*, \alpha_{J_s}\big)}=\sqrt{-1}\epsilon_{s_k}\mu_{s_k}l_{s_k}
\end{equation}
and for $s_k\in J_s$ 
\begin{equation}\label{3.4}
\frac{\partial U}{\partial \alpha_{s_k}}\bigg|_{z_s^*}
=\frac{\partial W}{\partial \alpha_{s_k}}\bigg|_{\big(a_{I_s}^*, \alpha_{J_s}\big)}=-\sqrt{-1}\mu_{s_k}l_{s_k}.
\end{equation}
By (\ref{3.3}) and (\ref{3.4}), we have  
\begin{equation}\label{3.26}
\begin{split}
\sum_{s=1}^c\bigg(\sum_{k=1}^6\frac{\partial U}{\partial \alpha_{s_k}}\bigg|_{z_s^*}\bigg)=&\sqrt{-1}\sum_{i\in I}\sum_{s_k\sim i} \epsilon_{s_k}\mu_{s_k}l_{s_k}-\sqrt{-1}\sum_{j\in J}\sum_{s_k\sim j} \mu_{s_k}l_{s_k}\\
=&\sqrt{-1} \sum_{i\in I}\epsilon_i\mu_i\Big(\sum_{s_k\sim i}l_{s_k}\Big)-\sqrt{-1}\sum_{j\in J}\mu_j\Big(\sum_{s_k\sim j} l_{s_k}\Big)\\
=&\sqrt{-1} \sum_{i\in I}\epsilon_i\mu_il_i-\sqrt{-1}\sum_{j\in J}\mu_jl_j\\
=&-\sqrt{-1}\sum_{i\in I}\epsilon_i\mu_i\Big(\mathrm H(v_i)-\frac{\iota_i}{2}\mathrm H(u_i)\Big)-\sqrt{-1}\sum_{j\in J}\mu_jl_j,
\end{split}
\end{equation}
where $s_k\sim i$ if $e_{s_k}$ intersects $L_i$ for $i\in I$ and $s_k\sim j$ if $e_{s_k}$ intersects $L_j$ for $j\in J,$ and the last equality comes from that \eqref{Sect3HuHv}.

Next, recall that for each $i\in I,$ $\alpha^*_i=\pi+\frac{\epsilon_i\mu_i\sqrt{-1}\mathrm H(u_i)}{2}$ and $\beta_i=\pi+\frac{\mu_i\theta_i}{2},$ and for each $j\in J,$ $\alpha_j=\pi+\frac{\mu_j\theta_j}{2}.$ \
Then for  $i\in I,$ we have
\begin{equation}\label{3.27}
\begin{split}
&\sqrt{-1} \bigg(\Big(p_i+\frac{\iota_i}{2}\Big)\alpha^*_i +\epsilon_i\beta_i \bigg)-\frac{\epsilon_i\mu_i}{2}\Big(\mathrm H(v_i)-\frac{\iota_i}{2}\mathrm H(u_i)\Big)\\
=&\sqrt{-1} \bigg(\Big(p_i+\frac{\iota_i}{2}\Big)\Big(\pi+\frac{\epsilon_i\mu_i\sqrt{-1}\mathrm H(u_i)}{2}\Big)+\epsilon_i\Big(\pi+\frac{\mu_i\theta_i}{2}\Big)\bigg)-\frac{\epsilon_i\mu_i}{2}\Big(\mathrm H(v_i)-\frac{\iota_i}{2}\mathrm H(u_i)\Big)\\
=&\Big(p_i+\frac{\iota_i}{2}+\epsilon_i\Big)\sqrt{-1}\pi+\frac{\epsilon_i\mu_i}{2}\Big(-p_i\mathrm H(u_i)-\frac{\iota_i}{2}\mathrm H(u_i)+\sqrt{-1}\theta_i-\mathrm H(v_i)+\frac{\iota_i}{2}\mathrm H(u_i)\Big)\\
=&\Big(p_i+\frac{\iota_i}{2}+\epsilon_i\Big)\sqrt{-1}\pi,
\end{split}
\end{equation}
where the last equality comes from $p_i\mathrm H(u_i)+\mathrm H(v_i)=\sqrt{-1}\theta_i.$
For $i\in I,$  we also have
\begin{equation}\label{3.28}
\begin{split}
\sqrt{-1}\big(q_i\beta_i+\epsilon_i\alpha_i^*\big)=&\sqrt{-1}\bigg(q_i\Big(\pi+\frac{\mu_i\theta_i}{2}\Big)+\epsilon_i\Big(\pi+\frac{\epsilon_i\mu_i\sqrt{-1}\mathrm H(u_i)}{2}\Big)\bigg)\\
=&\big(q_i+\epsilon_i\big)\sqrt{-1}\pi+\frac{\mu_i}{2}\Big(q_i\sqrt{-1}\theta_i-\mathrm H(u_i)\Big)\\
=&\big(q_i+\epsilon_i\big)\sqrt{-1}\pi+\frac{\mu_i}{2}\Big(q_i\big(p_i\mathrm H(u_i)+\mathrm H(v_i)\big)-\mathrm H(u_i)\Big)\\
=&\big(q_i+\epsilon_i\big)\sqrt{-1}\pi+\frac{\mu_i}{2}\mathrm H(\gamma_i).
\end{split}
\end{equation}
For each $j\in J,$ we have
\begin{equation}\label{3.29}
\begin{split}
\sqrt{-1}\Big(p_j+\frac{\iota_j}{2}\Big)\alpha_j-\frac{\mu_j}{2}l_j=&\sqrt{-1}\Big(p_j+\frac{\iota_j}{2}\Big)\Big(\pi+\frac{\mu_j\theta_j}{2}\Big)-\frac{\mu_j}{2}l_j\\
=&\Big(p_j+\frac{\iota_j}{2}\Big)\sqrt{-1}\pi+\frac{\mu_j}{2}\Big(p_j\sqrt{-1}\theta_j+\frac{\iota_j}{2}\sqrt{-1}\theta_j-l_j\Big)\\
=&\Big(p_j+\frac{\iota_j}{2}\Big)\sqrt{-1}\pi+\frac{\mu_j}{2}\Big(p_j\mathrm H(u_j)+\mathrm H(v_j)\Big)\\
=&\Big(p_j+\frac{\iota_j}{2}\Big)\sqrt{-1}\pi+\frac{\mu_j}{2}\mathrm H(\gamma_j).
\end{split}
\end{equation}
Then  for the case that $M_c$ is with a hyperbolic cone metric with singular locus $L_{\text{FSL}},$  the result follows from (\ref{L3.3}), (\ref{3.27}), (\ref{3.28}) and (\ref{3.29}).
\end{proof}

\begin{lemma}\label{L3.5} For $s\in \{1,\dots, c\},$ let $u_{s_1},\dots, u_{s_6}$ be the meridians of a tubular neighborhood of the components of $L_{\text{FSL}}$ intersecting the six edges of $\Delta_s.$ Then
\begin{equation}\label{CM0} 
\begin{split}
&\frac{e^{-{\sqrt{-1}}\sum_{s_i\in I_s}\alpha^*_{s_i}-{\sqrt{-1}}\sum_{s_j\in J_s}\alpha_{s_j}+4\sqrt{-1}\xi_s^*-\sum_{i=1}^4\log\big(1-e^{2\sqrt{-1}(\xi^*_s-\tau^*_{s_i})}\big)}}{{\frac{\partial ^2U}{\partial \xi_s^2}\Big|_{z^*_s}}}\\&\quad\quad\quad\quad\quad\quad\quad\quad\quad\quad\quad\quad\quad\quad\quad\quad\quad\quad\quad\quad=\frac{-1}{16\sqrt{\det\mathbb G\bigg(\frac{\mathrm H(u_{s_1})}{2},\dots,\frac{\mathrm H(u_{s_6})}{2}\bigg)}}.
\end{split}
\end{equation}
\end{lemma}

\begin{proof} The proof follows the same argument as in the proof of \cite[Lemma 3]{CM}. For $s_i\in I_s$ let $u_i=e^{\sqrt{-1}\alpha^*_{s_i}},$ and for $s_j\in J_s$ let $u_j=e^{\sqrt{-1}\alpha_{s_j}}.$ Let $z=e^{-2\sqrt{-1}\xi^*_s}$ and let $z'$ be the other root of equation (\ref{qe}). Then by (\ref{partialUxi}),
\begin{equation*}
\begin{split}
\frac{\partial ^2U}{\partial \xi_s^2}\bigg|_{z^*_s}=-4\bigg(&\frac{z}{1-z}+\frac{zu_1u_2u_4u_5}{1-zu_1u_2u_4u_5}+\frac{zu_1u_3u_4u_6}{1-zu_1u_3u_4u_6}+\frac{zu_2u_3u_5u_6}{1-zu_2u_3u_5u_6}\\
& -\frac{zu_1u_2u_3}{1-zu_1u_2u_3}-\frac{zu_1u_5u_6}{1-zu_1u_5u_6}-\frac{zu_2u_4u_6}{1-zu_2u_4u_6}-\frac{zu_3u_4u_5}{1-zu_3u_4u_5}\bigg).
\end{split}
\end{equation*}
Let $\mathrm{LHS}$ be the left hand side of (\ref{CM0}). Then we have
\begin{equation*}
\begin{split}
\frac{1}{\mathrm{LHS}}=&\frac{-4(1-zu_1u_2u_3)(1-zu_1u_5u_6)(1-zu_2u_4u_6)(1-zu_3u_4u_5)}{z^2u_1u_2u_3u_4u_5u_6}\\
&\times\bigg(\frac{z}{1-z}+\frac{zu_1u_2u_4u_5}{1-zu_1u_2u_4u_5}+\frac{zu_1u_3u_4u_6}{1-zu_1u_3u_4u_6}+\frac{zu_2u_3u_5u_6}{1-zu_2u_3u_5u_6}\\
&\quad-\frac{zu_1u_2u_3}{1-zu_1u_2u_3}-\frac{zu_1u_5u_6}{1-zu_1u_5u_6}-\frac{zu_2u_4u_6}{1-zu_2u_4u_6}-\frac{zu_3u_4u_5}{1-zu_3u_4u_5}\bigg).
\end{split}
\end{equation*}
By a direction computation and (\ref{ece}),  (\ref{qe}), (\ref{za}), (\ref{za'}) and (\ref{detG}), we have
\begin{equation*}
\begin{split}
\frac{1}{\mathrm{LHS}}=&-4\bigg(3Az+2B+\frac{C}{z}\bigg)=-4\bigg(Az-\frac{C}{z}\bigg)\\
=&-4A(z-z')=-4\sqrt{B^2-4AC}\\
=&-16\sqrt{\det\mathbb G\bigg(\Big(\sqrt{-1}(\alpha^*_{s_i}-\pi)\Big)_{s_i\in I_s},\Big(\sqrt{-1}(\alpha_{s_j}-\pi)\Big)_{s_j\in J_s}\bigg)}\\
=&-16\sqrt{\det\mathbb G\bigg(\frac{\mathrm H(u_{s_1})}{2},\dots,\frac{\mathrm H(u_{s_6})}{2}\bigg)},
\end{split}
\end{equation*}
from which (\ref{CM}) follows.
\end{proof}

\begin{proof}[Proof of Theorem \ref{main1'}]  
By (\ref{CRT}), Lemmas \ref{L3.1}, \ref{L3.2}, \ref{L3.4}, \ref{L3.5} and Theorem  \ref{Tor1} (2), we have
\begin{equation*}
\begin{split}
\frac{C^{{ \epsilon_I}}}{\sqrt{-\det\mathrm{Hess}\frac{\mathcal W^{{ \epsilon_I}}(z^{{ \epsilon_I}})}{4\pi\sqrt{-1}}}}=&\frac{(-1)^{\sum_{i\in I} \big(p_i+\frac{\iota_i}{2}+q_i\big)+\sum_{j\in J}\big(p_j+\frac{\iota_j}{2}\big)-\frac{rc}{2}}r^{\frac{|I|-c}{2}}}{2^{\frac{3|I|-c}{2}}\pi^{\frac{|I|+c}{2}}\big(4\pi\sqrt{-1}\big)^{-\frac{|I|+c}{2}}} \\
&\frac{e^{\frac{1}{2}\sum_{k=1}^n\mu_k\mathrm H(\gamma_k)}}{\sqrt{-(-16)^{c} (-2)^{|I|}\det\bigg(\frac{\partial \mathrm H(m_{i_1})}{\partial \mathrm H(u_{i_2})}\bigg)_{i_1,i_2\in I} \prod_{s=1}^{c}\sqrt{\det\mathbb G\bigg(\frac{\mathrm H(u_{s_1})}{2},\dots,\frac{\mathrm H(u_{s_6})}{2}\bigg)}}}\\
=&\frac{(-1)^{\sum_{i\in I} \big(p_i+\frac{\iota_i}{2}+q_i\big)+\sum_{j\in J}\big(p_j+\frac{\iota_j}{2}\big)-\frac{rc}{2}-\frac{|I|+c}{4}}}{2^{|I|-c}} \frac{e^{\frac{1}{2}\sum_{k=1}^n\mu_k\mathrm H(\gamma_k)}}{\sqrt{\pm\mathbb{T}_{(M\setminus L,\mathbf m)}([\rho_{M_{L_\theta}}])}}  r^{\frac{|I|-c}{2}},
\end{split}
\end{equation*}
where $\mathbf m=\big((m_i)_{i\in I}, (u_j)_{j\in J}\big).$
Therefore,
\begin{equation*}
\begin{split}
&\kappa_r\bigg(\sum_{\epsilon_I\in\{1,-1\}^{I}}\frac{C^{{ \epsilon_I}}}{\sqrt{-\det\mathrm{Hess}\frac{\mathcal W^{{ \epsilon_I}}(z^{{ \epsilon_I}})}{4\pi\sqrt{-1}}}}\bigg)=C \frac{e^{\frac{1}{2}\sum_{k=1}^n\mu_k\mathrm H(\gamma_k)}}{\sqrt{\pm\mathbb{T}_{(M\setminus L,\mathbf m)}([\rho_{M_{L_\theta}}])}},
\end{split}
\end{equation*}
where 
$$C={(-1)^{\sum_{i\in I} \big(p_i+\frac{\iota_i}{2}+q_i\big)+\sum_{j\in J}\big(p_j+\frac{\iota_j}{2}\big)-\sigma(L'\cup L_I)\big(-\frac{3}{r}-\frac{r+1}{4}\big)-\frac{r}{4}\big(\sum_{i\in I}q_i+\sum_{k=1}^np_k\big)+\frac{|I|}{2}+(|I|-c)\big(\frac{r}{2}-\frac{1}{4}\big)}}$$
is a quantity of norm $1,$ which is independent of the geometric structure on $M.$ 

Finally, by Theorem \ref{rrt}, we have 
$$\mathrm{RT}_r(M,L,(\mathbf b^{(r)}_I,\mathbf a^{(r)}_J))=C \frac{e^{\frac{1}{2}\sum_{k=1}^n\mu_k\mathrm H(\gamma_k)}}{\sqrt{\pm\mathbb{T}_{(M\setminus L,\mathbf m)}([\rho_{M_{L_\theta}}])}}e^{\frac{r}{4\pi}\big(\mathrm{Vol}(M_{L_\theta})+\sqrt{-1}\mathrm{CS}(M_{L_\theta})\big)}\Big( 1 + O \Big( \frac{1}{r} \Big) \Big).$$
\end{proof}


\section{Asymptotic expansion of the relative Turaev-Viro invariants}

In \cite{Y},  the second author studied the exponential growth rate of the relative Turaev-Viro invariants of an ideally triangulated $3$-manifold $N$ with boundary, and related it to the volume of the hyperbolic  metric on $N$ with the cone angles determined by the sequence of the colorings. The proof of Theorem \ref{main2} is based on the results in \cite{Y}. Throughout this section, we follow the conventions and normalizations of the relative Turaev-Viro invariants in \cite{Y}.

\subsection{Hyperbolic polyhedral $3$-manifolds}

Let $N$ be a compact $3$-manifold with non-empty boundary, and let $\mathcal T$ be an ideal triangulation of $N,$ that is, a finite collection  $T=\{\Delta_1,\dots,\Delta_{|T|}\}$ of truncated Euclidean tetrahedra with faces identified in pairs by affine homeomorphisms. We also let $E=\{e_1,\dots,e_{|E|}\}$ be the set of edges of $\mathcal T.$ As defined in \cite{Luo, LY}, a \emph{hyperbolic polyhedral metric} on  $(N,\mathcal T)$ is obtained by replacing each tetrahedron in $\mathcal T$ by a truncated hyperideal tetrahedron and replacing the gluing homeomorphisms between pairs of the faces  by isometries. The \emph{cone angle} at an edge is the sum of the dihedral angles of the truncated hyperideal tetrahedra around the edge. If all the cone angles are equal to $2\pi,$ then the hyperbolic polyhedral metric gives a hyperbolic metric on $M$ with totally geodesic boundary.  In \cite[Theorem 1.2 (b)]{LY}, Luo and the second author proved that hyperbolic polyhedral metrics on $(N, \mathcal T )$ are rigid in the sense that they are up to isometry determined by their cone angles.


\subsection{Growth rate of the relative Turaev-Viro invariants}
In this section, we recall the results from \cite{Y} on the exponential growth rate of the relative Turaev-Viro invariants of an ideally triangulated $3$-manifold. Suppose $N$ is a $3$-manifold with non-empty boundary and $\mathcal T$ is an ideal triangulation of $N$ with the set of edges $E$ and the set of tetrahedra $T.$ For a positive integer $r\geqslant 3,$ a \emph{coloring} $\mathbf a$ of $(N,\mathcal T)$ assigns an integer $a_i$ in between $0$ and $r-2$ to the edge $e_i,$
Let $\{\mathbf b^{(r)}\}$ be a sequence of colorings of $(N,\mathcal T)$ such that for each $i\in \{1,\dots, |E|\},$ either $b_i^{(r)} > \frac{r}{2}$ for all $r$ or  $b_i^{(r)} < \frac{r}{2}$ for all $r.$ In the former case we let $\mu_i=1$ and in the latter case we let $\mu_i=-1,$ and we let 
$$\theta_i^{(r)}=\mu_i\bigg(\frac{4\pi b_i^{(r)}}{r}-2\pi\bigg)=\bigg| \frac{4\pi b_i^{(r)}}{r}-2\pi\bigg|.$$
For simplicity, in the rest of this section we will drop the superscript  and write $\theta_i=\theta_i^{(r)},$ and keep in mind that it is a quantity depending on $r.$  
Let $\theta=(\theta_1,\dots,\theta_{|E|}),$ and let $N_{E_\theta}$ be $N$ with the hyperbolic polyhedral metric with cone angles $\theta,$ if exists, and let $\mathrm{Vol}(N_{E_\theta})$ be the volume of $N_{E_\theta}.$ 

\begin{theorem}(\cite[Theorem 1.4, Propositions 4.2, 4.3, 5.2, 5.5, 5.9 and 5.10]{Y})\label{rtv} Let $\mathrm {TV}_r(N,E,\mathbf b^{(r)})$ be the $r$-th relative Turaev-Viro invariant of $(N,\mathcal T)$ with the coloring $\mathbf b^{(r)}$ on edges. Suppose $ \theta_1,\dots,\theta_{|E|}$ converge to sufficiently small limits. Then as $r$ varies over all positive odd integers and at $q=e^{\frac{2\pi\sqrt{-1}}{r}},$ 
\begin{equation}\label{RTV}
\begin{split}
\mathrm {TV}_r(N,E,\mathbf b^{(r)})=&\Bigg( \sum_{\epsilon\in\{1,-1\}^{E}}\frac{C^{\epsilon}}{\sqrt{-\det\mathrm{Hess}\frac{\mathcal W^{\epsilon}(z^{\epsilon})}{4\pi\sqrt{-1}}}}\Bigg)  e^{\frac{r}{2\pi}\mathrm{Vol}(N_{E_\theta})}\Big( 1 + O \Big( \frac{1}{r} \Big) \Big),
\end{split}
\end{equation}
where $\mathcal W^\epsilon,$ $z^\epsilon$ and $C^\epsilon$ are quantities depending on $r$ which will be explained as follows.
\end{theorem}

For the function $\mathcal W^\epsilon,$ we let $\alpha_i=\frac{2\pi a_i}{r}$ and $\beta_i=\frac{2\pi b^{(r)}_i}{r}$ for $i\in\{1,\dots,|E|\},$ $\xi_s=\frac{2\pi k_s}{r}$ for $s\in\{1,\dots, |T|\},$ $\tau_{s_i}=\frac{2\pi T_{s_i}}{r}$ for $i\in\{1,\dots,4\},$ and $\eta_{s_j}=\frac{2\pi Q_{s_j}}{r}$ for $j\in\{1,2,3\}.$ Let
$$\mathrm {D_A}=\Big\{(\alpha,\xi)\in\mathbb R^{|E|+|T|}\ \Big|\ (\alpha_{s_1},\dots,\alpha_{s_6}) \text{ is admissible, } \max\{\tau_{s_i}\}\leqslant \xi_s\leqslant \min\{\eta_{s_j}, 2\pi\}, s\in\{1,\dots, |T|\}\Big\}$$
where $\alpha=(\alpha_1,\dots,\alpha_{|E|}),$ and let
$$\mathrm {D_H}=\Big\{(\alpha,\xi)\in\mathrm {D_A} \ \Big|\ (\alpha_{s_1},\dots,\alpha_{s_6}) \text{ is of the hyperideal type}, s\in\{1,\dots, |T|\} \Big\}.$$
For a sufficiently small $\delta >0,$ let 
$$\mathrm {D_H^\delta}=\Big\{(\alpha,\xi)\in\mathrm {D_H}\ \Big|\ d((\alpha,\xi), \partial\mathrm {D_H})>\delta \Big\},$$
where $d$ is the Euclidean distance on $\mathbb R^n.$
Then $\mathcal W^\epsilon$ is the following function 
\begin{equation}\label{WforTV}
\mathcal W^{\epsilon}(\alpha,\xi)=-\sum_{i=1}^{|E|}2\epsilon_i(\alpha_i-\pi)(\beta_i-\pi)+\sum_{s=1}^{|T|}U(\alpha_{s_1},\dots,\alpha_{s_6},\xi_s)
\end{equation}
where $U$ is as defined in (\ref{term}), which  is continuous on 
$$\mathrm{D_{H,\mathbb C}}=\big\{(\alpha,\xi)\in\mathbb C^{|E|+|T|}\ \big|\ (\mathrm{Re}(\alpha),\mathrm{Re}(\xi))\in \mathrm{D_{H}}\big\}$$ and for any $\delta>0$ is analytic on 
$$\mathrm{D^\delta_{H,\mathbb C}}=\big\{(\alpha,\xi)\in\mathbb C^{|E|+|T|}\ \big|\ (\mathrm{Re}(\alpha),\mathrm{Re}(\xi))\in \mathrm{D^\delta_{H}}\big\},$$ 
where $\mathrm{Re}(\alpha)=(\mathrm{Re}(\alpha_1),\dots, \mathrm{Re}(\alpha_{|E|}))$ and $\mathrm{Re}(\xi)=(\mathrm{Re}(\xi_1),\dots, \mathrm{Re}(\xi_{|T|})).$

For $z^\epsilon,$ we for each $i\in \{1,\dots,|E|\}$ let $l_i$ be the length of the edge $e_i$ in $N_{E_\theta}$ and let
\begin{align}\label{defnalphaitv}
\alpha^*_i=\pi+\epsilon_i\mu_i\sqrt{-1}l_i.    
\end{align}
For each $s\in\{1,\dots,|T|\},$  let 
$$\xi^*_s=\xi(\alpha^*_{s_1},\dots,\alpha^*_{s_6}).$$
 Suppose $ \theta_1,\dots,\theta_{|E|}$ are sufficiently small.  Then by \cite[Proposition 5.2]{Y}, $\mathcal W^{\epsilon}$ has a unique critical point 
$$z^{\epsilon}=\big(\alpha^*_1,\dots,\alpha^*_{|E|},\xi^*_1,\dots,\xi^*_{|T|}\big)$$
in $\mathrm{D_{H,\mathbb C}}$ with critical value 
$$2|T|\pi^2+2\sqrt{-1}\mathrm{Vol}(N_{E_\theta}).$$

By the computation in the end of the proof of \cite[Proposition 5.5]{Y}, we have
\begin{equation}\label{CTV}
C^\epsilon=2^{\mathrm{rank H}_2(N;\mathbb Z_2)}\frac{(-1)^{|E|+\frac{r}{2}(|E|-|T|)}r^{\frac{|E|-|T|}{2}}}{2^{\frac{3|E|+|T|}{2}}\pi^{\frac{|E|+|T|}{2}}\{1\}^{|E|-{|T|}}}e^{\sum_{i=1}^{|E|}\epsilon_i\sqrt{-1}(\alpha^*_i+\beta_i)+\sum_{s=1}^{|T|}\kappa(\alpha_s^*,\xi_s^*)},
\end{equation}
where $\kappa$ is as defined in (\ref{kappa}).

Finally, we notice that all of $\mathcal W^\epsilon,$ $z^\epsilon$ and $C^\epsilon$ depend on $r$ because each $\beta_i=\frac{2\pi b_i^{(r)}}{r},$ $i\in\{1,\dots,|E|\},$ does. 


\subsection{Proof of Theorem \ref{main2}}

Note that the function $\mathcal W^{\epsilon}(\alpha,\xi)$ in \eqref{WforTV} equals the function $\mathcal W^{\epsilon_I}(\alpha_I,\xi)$ in \eqref{WforRT} with $q_i=\iota_i=p_j=0$, $I=E$, $J=\emptyset$ and $c=|T|$. 
By Theorem \ref{rtv}, to prove Theorem \ref{main2}, it suffices to compute  $ \frac{C^{\epsilon}}{\sqrt{-\det\mathrm{Hess}\frac{\mathcal W^{\epsilon}(z^{\epsilon})}{4\pi\sqrt{-1}}}}$ for each $\epsilon\in\{-1,1\}^E.$ To do this, we need the following Lemmas  \ref{L4.2}, \ref{L4.3} and \ref{L4.4}, which are analogues of Lemma \ref{L3.1}, \ref{L3.2} and \ref{L3.5} respectively.

\begin{lemma}\label{L4.2} For $s\in\{1,\dots,|T|\},$ let $z^*_s=(\alpha^*_{s_1},\dots,\alpha^*_{s_6},\xi^*_s).$ Then
$$-\det\mathrm{Hess}\mathcal W^{\epsilon}(z^{\epsilon})=-(-1)^{\frac{3|E|}{2}}\det\bigg(\frac{\partial \theta_i}{\partial l_j}\bigg)_{ij} \prod _{s=1}^{|T|} \frac{\partial ^2U}{\partial \xi_s^2}\bigg|_{z^*_s}.$$
\end{lemma}

\begin{proof} 
For $s\in\{1,\dots,|T|\}$ and $i\in \{1,\dots,|E|\},$ we denote by $s\sim i$ if the tetrahedron $\Delta_s$ intersects the edge $e_i;$  and for $\{i,j\}\subset\{1,\dots,|E|\}$ we denote by $s\sim i,j$ if $\Delta_s$ intersects both $e_i$ and $e_j.$ For $s\in\{1,\dots,|T|\},$ let $\alpha_s=(\alpha_{s_1},\dots,\alpha_{s_6})$ and let $\alpha_s^*=(\alpha^*_{s_1},\dots,\alpha^*_{s_6}).$ 
We claim that
\begin{enumerate}[(1)]
\item For $s\in\{1,\dots,|T|\},$
$$\frac{\partial^2 \mathcal W^\epsilon}{\partial \xi_s^2}\bigg|_{z^\epsilon}=\frac{\partial^2 U}{\partial \xi_s^2}\bigg|_{z_s^*}.$$

\item For $\{s_1,s_2\}\subset\{1,\dots,|T|\},$
$$\frac{\partial^2 \mathcal W^\epsilon}{\partial \xi_{s_1}\partial \xi_{s_2}}\bigg|_{z^\epsilon}=0.$$

\item For $i\in\{1,\dots,|E|\}$ and $s\in\{1,\dots,|T|\},$
$$\frac{\partial^2 \mathcal W^\epsilon}{\partial \alpha_i\partial \xi_s}\bigg|_{z^\epsilon}=-\frac{\partial^2 U}{\partial \xi_s^2}\bigg|_{z^*_s} \frac{\partial \xi_s(\alpha_s)}{\partial \alpha_i}\bigg|_{\alpha_s^*}.$$
\item For $i\in\{1,\dots,|E|\},$ 
$$\frac{\partial^2 \mathcal W^\epsilon}{\partial \alpha_i^2}\bigg|_{z^\epsilon}=-\sqrt{-1}\frac{\partial \theta_i}{\partial l_i}+\sum_{s\sim i}\frac{\partial^2 U}{\partial \xi_s^2}\bigg|_{z^*_s} \bigg(\frac{\partial \xi_s(\alpha_s)}{\partial \alpha_i}\bigg|_{\alpha_s^*}\bigg)^2.$$

\item[(5)] For $\{i,j\}\subset\{1,\dots,|E|\},$ 
$$\frac{\partial^2 \mathcal W^\epsilon}{\partial \alpha_i\partial \alpha_j}\bigg|_{z^\epsilon}=-\sqrt{-1}\frac{\epsilon_i\mu_i}{\epsilon_j\mu_j} \frac{\partial \theta_i}{\partial l_j}+\sum_{s\sim i,j}\frac{\partial^2 U}{\partial \xi_s^2}\bigg|_{z^*_s} \frac{\partial \xi_s(\alpha_s)}{\partial \alpha_i}\bigg|_{\alpha_s^*} \frac{\partial \xi_s(\alpha_s)}{\partial \alpha_j}\bigg|_{\alpha_s^*}.$$
\end{enumerate}
Assuming (1)-(5),  the result follows from statements analogous to (\ref{congurent0}), (\ref{detD0}) and (\ref{detA0}).
The proofs of (1)-(3) follow from the same argument as in the proof of Lemma \ref{L3.1}. For (4) and (5), we have 
\begin{equation}\label{4.10}
\frac{\partial ^2\mathcal W^\epsilon}{\partial \alpha_i^2}\bigg|_{z^\epsilon}=\sum_{s\sim i}\frac{\partial ^2U}{\partial \alpha_i^2}\bigg|_{z^*_s},
\end{equation}
and
\begin{equation}\label{4.11}
\frac{\partial ^2\mathcal W^\epsilon}{\partial \alpha_i\partial \alpha_j}\bigg|_{z^\epsilon}=\sum_{s\sim i,j}\frac{\partial ^2U}{\partial \alpha_i\partial \alpha_j}\bigg|_{z^*_s}.
\end{equation}
From \eqref{3.13}, for $i,j\in\{1,\dots, 6\}$ we have
\begin{equation}\label{4.13}
\begin{split}
\frac{\partial ^2U}{\partial \alpha_{s_i}\partial \alpha_{s_j}}\bigg|_{(\alpha_s,\xi_s(\alpha_s))}
=&\frac{\partial^2 W}{\partial \alpha_{s_i}\partial \alpha_{s_j}}\bigg|_{\alpha_s}+\frac{\partial^2 U}{\partial \xi_s^2}\bigg|_{(\alpha_s,\xi_s(\alpha_s))} \frac{\partial \xi_s(\alpha_s)}{\partial \alpha_{s_i}}\bigg|_{\alpha_s}  \frac{\partial \xi_s(\alpha_s)}{\partial \alpha_{s_j}}\bigg|_{\alpha_s}.
\end{split}
\end{equation}
By Theorem \ref{co-vol}, we have
\begin{equation*}
\frac{\partial^2 W}{\partial \alpha_{s_i}\partial \alpha_{s_j}}\bigg|_{\alpha^*_s}=-\sqrt{-1}\frac{\epsilon_{s_i}\mu_{s_i}}{\epsilon_{s_j}\mu_{s_j}} \frac{\partial \theta_{s_i}}{\partial l_{s_j}},
\end{equation*}
where $\theta_{s_k}$ is the dihedral angle of  $\Delta_s$ at $e_{s_k},$ and $\epsilon_{s_k}=\epsilon_i,$ $\mu_{s_k}=\mu_i$ and $l_{s_k}=l_i$ if the edge $e_{s_k}$ intersections the edge $e_i.$ Then (4) and (5) follow from (\ref{4.10}),  (\ref{4.11}), (\ref{4.13}) and the fact that
$$\theta_i=\sum_{s\sim i}\sum_{s_k}\theta_{s_k},$$
where the second sum is over $s_k$ such that the edge $e_{s_k}$ in $\Delta_s$ intersects the edge $e_i.$ 
\end{proof}

\begin{lemma}\label{L4.3} For $s\in\{1,\dots,|T|\},$ let $\alpha_s^*=(\alpha^*_{s_1},\dots,\alpha^*_{s_6}).$ Then
\begin{equation*}
\begin{split}
\kappa(\alpha_s^*,\xi_s^*)=&-\frac{\sqrt{-1}}{2}\sum_{i=1}^6\epsilon_{s_i}\mu_{s_i}\theta_{s_i}\\
&-\frac{\sqrt{-1}}{2}\sum_{i=1}^6\alpha^*_{s_i}+2\sqrt{-1}\xi_s^*-\frac{1}{2}\sum_{i=1}^4\log\big(1-e^{2\sqrt{-1}(\xi^*_s-\tau^*_{s_i})}\big),
\end{split}
\end{equation*}
where
$\tau_{s_1}=\frac{\alpha^*_{s_1}+\alpha^*_{s_2}+\alpha^*_{s_3}}{2},$ $\tau_{s_2}=\frac{\alpha^*_{s_1}+\alpha^*_{s_5}+\alpha^*_{s_6}}{2},$ $\tau_{s_3}=\frac{\alpha^*_{s_2}+\alpha^*_{s_4}+\alpha^*_{s_6}}{2}$ and $\tau_{s_4}=\frac{\alpha^*_{s_3}+\alpha^*_{s_4}+\alpha^*_{s_5}}{2}.$ 
\end{lemma}

\begin{proof} For $s\in\{1,\dots,|T|\},$ let $z^*_s=(\alpha^*_{s_1},\dots,\alpha^*_{s_6},\xi^*_s).$ Then by (\ref{kappapf2}), we have
\begin{equation*}
\kappa(\alpha_s^*,\xi_s^*)=-\frac{\sqrt{-1}}{2}\sum_{i=1}^6\frac{\partial U}{\partial \alpha_{s_i}}\bigg|_{z^*_s}-\frac{\sqrt{-1}}{2}\sum_{i=1}^6\alpha^*_{s_i}+2\sqrt{-1}\xi_s^*-\frac{1}{2}\sum_{i=1}^4\log\big(1-e^{2\sqrt{-1}(\xi^*_s-\tau^*_{s_i})}\big).
\end{equation*}
For $i\in\{1,\dots,6\}$, by \eqref{3.5} and Theorem \ref{co-vol}, we have 
\begin{align*}
    \frac{\partial U}{\partial \alpha_{s_i}}\bigg|_{z^*_s}=
    \frac{\partial W}{\partial \alpha_{s_i}}\bigg|_{\alpha^*_s}=\epsilon_{s_i}\mu_{s_i}\theta_{s_i}.
\end{align*}
This implies the desired result.
\end{proof}

\begin{lemma}\label{L4.4}
\begin{equation}\label{CM}
\frac{e^{-{\sqrt{-1}}\sum_{i=1}^6\alpha^*_{s_i}+4\sqrt{-1}\xi_s^*-\sum_{i=1}^4\log\big(1-e^{2\sqrt{-1}(\xi^*_s-\tau^*_{s_i})}\big)}}{{\frac{\partial ^2U}{\partial \xi_s^2}\Big|_{z^*_s}}}=\frac{-1}{16\sqrt{\det\mathbb G(l_{s_1},\dots,l_{s_6})}}.
\end{equation}
\end{lemma}

\begin{proof} 
Following the computations in the proof of Lemma \ref{L3.5}, we have
\begin{equation*}
\begin{split}
\frac{1}{\mathrm{LHS}}
=&-16\sqrt{\det\mathbb G\bigg(\sqrt{-1}(\alpha^*_{s_1}-\pi),\dots,\sqrt{-1}(\alpha^*_{s_6}-\pi)\bigg)}\\
=&-16\sqrt{\det\mathbb G(l_{s_1},\dots,l_{s_6})},
\end{split}
\end{equation*}
where the last equality follows from \eqref{defnalphaitv}. This implies the desired result.
\end{proof}

\begin{proof}[Proof of Theorem \ref{main2}] We recall that for each $\epsilon\in\{1,-1\}^{E},$ $\alpha^*_i=\pi+\epsilon_i\mu_i\sqrt{-1}l_i$ and $\beta_i=\pi+\frac{\mu_i\theta_i}{2}.$ Then by (\ref{CTV}), Lemmas \ref{L4.2}, \ref{L4.3} and \ref{L4.4} and Theorem \ref{Tor2} (2), we have
\begin{equation*}
\begin{split}
\frac{C^{\epsilon}}{\sqrt{-\det\mathrm{Hess}\frac{\mathcal W^{\epsilon}(z^{\epsilon})}{4\pi\sqrt{-1}}}}=&2^{\mathrm{rank H}_2(N;\mathbb Z_2)}\frac{(-1)^{|E|+\frac{r}{2}(|E|-|T|)}r^{\frac{|E|-|T|}{2}}}{2^{\frac{3|E|+|T|}{2}}\pi^{\frac{|E|+|T|}{2}}\{1\}^{|E|-{|T|}}} \\
&\frac{e^{\sum_{i=1}^{|E|}\epsilon_i\sqrt{-1}\big(2\pi+\epsilon_i\mu_i\sqrt{-1}l_i+\frac{\mu_i\theta_i}{2}\big)+\sum_{s=1}^{|T|}\big(-\frac{\sqrt{-1}}{2}\sum_{i=1}^6\epsilon_{s_i}\mu_{s_i}\theta_{s_i}\big)}}{\big(4\pi\sqrt{-1}\big)^{-\frac{|E|+|T|}{2}}\sqrt{-(-1)^{\frac{3|E|}{2}}(-16)^{|T|} \det\Big(\frac{\partial \theta_i}{\partial l_j}\Big)_{ij} \prod_{s=1}^{|T|}\sqrt{\det\mathbb G(l_{s_1},\dots,l_{s_6})}}}\\
=&\frac{(-1)^{{\frac{3|E|+|T|}{4}+\frac{r}{2}(|E|-|T|)}}2^{\mathrm{rank H}_2(N;\mathbb Z_2)}}{2^{3|E|-2|T|}\pi^{|E|-|T|}} \frac{e^{-\sum_{i=1}^{|E|}\mu_il_i}}{\sqrt{\pm\mathbb{T}_{( M,\mathbf m)}([\rho_{M_{E_\theta}}])}}  r^{\frac{3|E|-3|T|}{2}}\Big( 1 + O\Big(\frac{1}{r}\Big)\Big),
\end{split}
\end{equation*}
where the last equality comes from the fact that
$$\sum_{i=1}^{|E|}\theta_i=\sum_{s=1}^{|T|}\sum_{i=1}^6\theta_{s_i}$$
and that
$$\{1\}=2\sqrt{-1}\sin\frac{2\pi}{r}=\frac{4\pi \sqrt{-1}}{r}+O\Big(\frac{1}{r^2}\Big).$$
Therefore, 
\begin{equation*}
\begin{split}
\sum_{\epsilon\in\{1,-1\}^{E}}&\frac{C^{\epsilon}}{\sqrt{-\det\mathrm{Hess}\frac{\mathcal W^{\epsilon}(z^{\epsilon})}{4\pi\sqrt{-1}}}}\\
=&\frac{(-1)^{{\frac{3|E|+|T|}{4}+\frac{r}{2}(|E|-|T|)}}2^{\mathrm{rank H}_2(N;\mathbb Z_2)}}{(4\pi)^{|E|-|T|}} \frac{e^{-\sum_{i=1}^{|E|}\mu_il_i}}{\sqrt{\pm\mathbb{T}_{( M,\mathbf m)}([\rho_{M_{E_\theta}}])}}  r^{\frac{3|E|-3|T|}{2}}\Big( 1 + O\Big(\frac{1}{r}\Big)\Big).
\end{split}
\end{equation*}

Finally, by Theorem \ref{rtv} and that $|E|-|T|=\chi(N),$ we have
\begin{equation*}
\begin{split}
\mathrm {TV}_r&(N,E,\mathbf b^{(r)})\\
=&\frac{(-1)^{|E|+\chi(N)\big(\frac{r}{2}-\frac{1}{4}\big)}2^{\mathrm{rank H}_2(N;\mathbb Z_2)}}{(4\pi)^{\chi(N)}}\frac{e^{-\sum_{i=1}^{|E|}\mu_il_i}}{\sqrt{\pm\mathbb{T}_{( M,\mathbf m)}([\rho_{M_{E_\theta}}])}}  r^{\frac{3}{2}\chi(N)}  e^{\frac{r}{2\pi}\mathrm{Vol}(N_{E_\theta})} \Big( 1 + O\Big(\frac{1}{r}\Big)\Big).
\end{split}
\end{equation*}
\end{proof}

\section{Asymptotic expansion of the discrete Fourier transforms of quantum $6j$-symbols}

In \cite{BY}, Belletti and the second author studied the exponential growth rate of the discrete Fourier transforms of the quantum $6j$-symbols and related it to the volume of the deeply truncated tetrahedron with dihedral angles determined by the sequence of the colorings. The proof of Theorem \ref{main3} is based on the results in \cite{BY}. Throughout this section, we follow the conventions and normalizations of the quantum $6j$-symbols and Yokota invariants in \cite{BY}.


\subsection{Growth rate of the discrete Fourier transforms of quantum $6j$-symbols}

In this section, we recall the results from \cite{BY} on the exponential growth rate of the exponential growth rate of the discrete Fourier transforms of the quantum $6j$-symbols. Let $(I,J)$ be a partition of $\{1,\dots,6\},$ and let $\Delta$ be a deeply truncated tetrahedron of type $(I,J),$ ie., $\{e_i\}_{i\in I}$ is the set of edges of deep truncation. For a $6$-tuple $(\mathbf b_I, \mathbf a_J)=((b_i)_{i\in I},(a_j)_{j\in J})$ of integers in $\{0,\dots, r-2\},$   let $\mathrm{\widehat {Y}}_r\big(\mathbf b_I; \mathbf a_J\big)$ be the discrete Fourier transform of the Yokota invariant of the trivalent graph $\cp{\includegraphics[width=0.5cm]{Yokota}}$ with respect to $(\mathbf b_I,\mathbf a_J).$

Let $\{(\mathbf b_I^{(r)}, \mathbf a_J^{(r)})\}$ be a sequence of $6$-tuples such that for any $i\in I,$ either $b_i^{(r)}>\frac{r}{2}$ for all $r$ or $b_i^{(r)}<\frac{r}{2}$ for all $r;$ and for any $j\in J,$ either $a_j^{(r)}>\frac{r}{2}$ for all $r$ or $a_j^{(r)}<\frac{r}{2}$ for all $r.$ In the former case we let $\mu_i=\mu_j=1$ and in the latter case we let $\mu_i=\mu_j=-1,$ and we let
$$\theta_i^{(r)}=\mu_i\bigg(\frac{2\pi b_i^{(r)}}{r}-\pi\bigg)=\bigg|\frac{2\pi b_i^{(r)}}{r}-\pi\bigg|$$
for $i\in I,$ and 
$$\theta_j^{(r)}=\mu_j\bigg(\frac{2\pi a_j^{(r)}}{r}-\pi\bigg)=\bigg|\frac{2\pi a_j^{(r)}}{r}-\pi\bigg|$$ 
for $j\in J.$ 
For simplicity, in the rest of this section we will drop the superscript  and write $\theta_i=\theta_i^{(r)}$ and $\theta_j=\theta_j^{(r)},$ and keep in mind that it is a quantity depending on $r.$  
Suppose $\Delta(\theta_I;\theta_J)$ is a deeply truncated tetrahedron  of type $(I,J)$ with $\theta_I=\{\theta_i\}_{i\in I}$ the set of dihedral angles at the edges of deep truncation and  $\theta_J=\{\theta_j\}_{j\in J}$ the set of dihedral angles at the regular edges, and $\mathrm{Vol}(\Delta(\theta_I;\theta_J))$ is the volume of $\Delta(\theta_I;\theta_J).$

\begin{theorem}(\cite[Theorem 1.2, Propositions 4.2, 4.3, 5.2, 5.5, 5.9 and 5.10]{BY})\label{dft}  Suppose $\theta_1,\dots,\theta_6$ converge to sufficiently small limits.
Then evaluated at the root of unity $q=e^{\frac{2\pi \sqrt{-1}}{r}}$ and as $r$ varies over all positive odd integers,
\begin{equation}\label{DF}
\begin{split}
\mathrm{\widehat {Y}}_r(\mathbf b^{(r)}_I; \mathbf a^{(r)}_J)=&\Bigg( \sum_{\epsilon_I\in\{1,-1\}^I}\frac{C^{\epsilon_I}}{\sqrt{-\det\mathrm{Hess}\frac{\mathcal W^{\epsilon_I}(z^{\epsilon_I})}{4\pi\sqrt{-1}}}}\Bigg)  e^{\frac{r}{\pi}\mathrm{Vol}(\Delta(\theta_I;\theta_J))}\Big( 1 + O \Big( \frac{1}{r} \Big) \Big),
\end{split}
\end{equation}
where $\mathcal W^{\epsilon_I},$ $z^{\epsilon_I}$ and $C^{\epsilon_I}$  are quantities depending on $r$ which will be explained as follows.
\end{theorem}

For the function $\mathcal W^{\epsilon_I},$ we let $\beta_i=\frac{2\pi b^{(r)}_i}{r}$ and $\alpha_i=\frac{2\pi a_i}{r}$ for $i\in I,$ $\alpha_i=\frac{2\pi a^{(r)}_j}{r}$ for $j\in J,$ $\xi_s=\frac{2\pi k_s}{r}$ for $s\in\{1,2\},$ $\tau_i=\frac{2\pi T_i}{r}$ for $i\in\{1,\dots,4\},$ and $\eta_j=\frac{2\pi Q_j}{r}$ for $j\in\{1,2,3\}.$ For a fixed $(\alpha_j)_{j\in J},$ let
$$\mathrm {D_A}=\Big\{(\alpha_I,\xi_1,\xi_2)\in\mathbb R^{|I|+2}\ \Big|\ (\alpha_1,\alpha_2,\dots,\alpha_6) \text{ is admissible, } \max\{\tau_i\}\leqslant \xi_s\leqslant \min\{\eta_j, 2\pi\}, s\in\{1,2\}\Big\}$$
where $\alpha_I=(\alpha_i)_{i\in I},$ and let
$$\mathrm {D_H}=\Big\{(\alpha_I,\xi_1,\xi_2)\in\mathrm {D_A} \ \Big|\ (\alpha_1,\alpha_2,\dots,\alpha_6) \text{ is of the hyperideal type} \Big\}.$$
For a sufficiently small $\delta >0,$ let 
$$\mathrm {D_H^\delta}=\Big\{(\alpha_I,\xi_1,\xi_2)\in\mathrm {D_H}\ \Big|\ d((\alpha_I,\xi_1,\xi_2), \partial\mathrm {D_H})>\delta \Big\},$$
where $d$ is the Euclidean distance on $\mathbb R^n.$ Then $\mathcal W^{\epsilon_I}$ is the following function
$$\mathcal W^{\epsilon_I}(\alpha_I,\xi_1,\xi_2)=-\sum_{i\in I}2\epsilon_i(\alpha_i-\pi)(\beta_i-\pi)+U(\alpha_1,\dots,\alpha_6,\xi_1)+U(\alpha_1,\dots,\alpha_6,\xi_2)$$
where $U$ is as defined in (\ref{term}), which is continuous on 
$$\mathrm{D_{H,\mathbb C}}=\big\{(\alpha_I,\xi_1,\xi_2)\in\mathbb C^{|I|+2}\ \big|\ (\mathrm{Re}(\alpha_I),\mathrm{Re}(\xi_1),\mathrm{Re}(\xi_2))\in \mathrm{D_{H}}\big\}$$ and for any $\delta>0$ is analytic on 
$$\mathrm{D^\delta_{H,\mathbb C}}=\big\{(\alpha_I,\xi_1,\xi_2)\in\mathbb C^{|I|+2}\ \big|\ (\mathrm{Re}(\alpha_I),\mathrm{Re}(\xi_1),\mathrm{Re}(\xi_2))\in \mathrm{D^\delta_{H}}\big\},$$ 
where $\mathrm{Re}(\alpha_I)=(\mathrm{Re}(\alpha_i))_{i\in I}.$

For $z^{\epsilon_I},$ we for each $i\in I$ let $l_i$ be the length of the edge of deep truncation $e_i$ and let
$$\alpha_i^*=\pi+\epsilon_i\mu_i\sqrt{-1}l_i;$$
and let 
$$\xi^*=\xi\big((\alpha_i^*)_{i\in I}, (\alpha_j)_{j\in J}\big).$$
 Suppose $\theta_1,\dots,\theta_6$ are sufficiently small. Then $\mathcal W^{\epsilon_I}(\alpha_I,\xi_1,\xi_2)$ has a critical point 
$$z^{\epsilon_I}=\Big((a_i^*)_{i\in I}, \xi^*, \xi^*\Big)$$
in $\mathrm{D_{H,\mathbb C}}$ with critical value 
$$4\pi^2+4\sqrt{-1}  \mathrm{Vol}\big(\Delta(\theta_I;\theta_J)\big).$$

By the computation in the end of the proof of \cite[Proposition 5.5]{BY}, we have
\begin{equation}\label{CDF}
\begin{split}
C^{\epsilon_I}=n(a_J)\frac{(-1)^{|I|+\frac{r}{2}(|I|-2)}r^{\frac{|I|-2}{2}}}{2^{\frac{3|I|+2}{2}}\pi^{\frac{|I|+2}{2}}\{1\}^{|I|-2}}e^{\sum_{i\in I}\epsilon_i\sqrt{-1}(\alpha^*_i+\beta_i)+2\kappa\big((\alpha_i^*)_{i\in I}, (\alpha_j)_{j\in J}, \xi^*\big)},
\end{split}
\end{equation}
where $\kappa$ is as defined in (\ref{kappa}), and $n(a_J)$ is the number of $3$-admissible colorings $c$ such that $c_j \equiv a_j\ (\text{mod } 2)$ for each $j\in J.$

Finally, we notice that all of $\mathcal W^\epsilon,$ $z^\epsilon$ and $C^\epsilon$ depend on $r$ because each $\beta_i=\frac{2\pi b_i^{(r)}}{r},$ $i\in I,$ and $\alpha_j=\frac{2\pi a_j^{(r)}}{r},$ $j\in J,$  does. 


\subsection{Proof of Theorem \ref{main3}}

By Theorem \ref{dft}, to prove Theorem \ref{main3}, it suffices to compute the $ \frac{C^{\epsilon_I}}{\sqrt{-\det\mathrm{Hess}\frac{\mathcal W^{\epsilon_I}(z^{\epsilon_I})}{4\pi\sqrt{-1}}}}$ for each $\epsilon_I\in \{1,-1\}^I.$ To do this, we need the following Lemmas  \ref{L5.2}, \ref{L5.3} and \ref{L5.4}
, which can be proved using the arguments from Lemmas \ref{L4.2}, \ref{L4.3} and \ref{L4.4} respectively. 

\begin{lemma}\label{L5.2} Let $a_I^*=(a_i^*)_{i\in I},$ $a_J=(a_j)_{j\in J}$ and $z^{\epsilon_I}=\big(a^*_I, \xi^*, \xi^*\big).$ Then
$$-\det\mathrm{Hess}\mathcal W^{\epsilon_I}(z^{\epsilon_I})=-(-1)^{\frac{3|I|}{2}}\det\bigg(\frac{\partial \theta_{i_1}}{\partial l_{i_2}}\bigg)_{i_1,i_2\in I} \frac{\partial ^2U}{\partial \xi_1^2}\bigg|_{\big(a^*_I, a_J, \xi^*\big)} \frac{\partial ^2U}{\partial \xi_2^2}\bigg|_{\big(a^*_I, a_J, \xi^*\big)}.$$
\end{lemma}

\begin{lemma}\label{L5.3} Let $\alpha_I^*=(\alpha^*_i)_{i\in I}$ and $\alpha_J=(\alpha_j)_{j\in J}.$ Then
\begin{equation*}
\begin{split}
\kappa(\alpha_I^*,\alpha_J,\xi^*)=&-\frac{\sqrt{-1}}{2}\sum_{i\in I}\epsilon_i\mu_i\theta_i-\frac{1}{2}\sum_{j\in J}\mu_jl_j\\
&-\frac{\sqrt{-1}}{2}\sum_{i\in I}\alpha^*_i-\frac{\sqrt{-1}}{2}\sum_{j\in J}\alpha_j+2\sqrt{-1}\xi^*-\frac{1}{2}\sum_{i=1}^4\log\big(1-e^{2\sqrt{-1}(\xi^*-\tau^*_i)}\big),
\end{split}
\end{equation*}
where with the notation that $\alpha^*_j=\alpha_j$ for $j\in J,$
$\tau_{1}=\frac{\alpha^*_{1}+\alpha^*_{ 2}+\alpha^*_{ 3}}{2},$ $\tau_{ 2}=\frac{\alpha^*_{ 1}+\alpha^*_{ 5}+\alpha^*_{ 6}}{2}$ and $\tau_{ 3}=\frac{\alpha^*_{ 2}+\alpha^*_{ 4}+\alpha^*_{ 6}}{2},$ $\tau_{ 4}=\frac{\alpha^*_{ 3}+\alpha^*_{ 4}+\alpha^*_{ 5}}{2}.$  \end{lemma}

\begin{lemma}\label{L5.4} Let $\alpha^*_I=(\alpha^*_i)_{i\in I},$ $l_I=(l_i)_{i\in I}$ and $\theta_J=(\theta_j)_{j\in J}.$ Then for $s\in\{1,2\},$
\begin{equation}\label{CM2} 
\frac{e^{-{\sqrt{-1}}\sum_{i\in I}\alpha^*_i-{\sqrt{-1}}\sum_{j\in J}\alpha_j+4\sqrt{-1}\xi^*-\sum_{i=1}^4\log\big(1-e^{2\sqrt{-1}(\xi^*-\tau^*_i)}\big)}}{{\frac{\partial ^2U}{\partial \xi_s^2}\Big|_{\big(\alpha_I^*, \xi^*\big)}}}=\frac{-1}{16\sqrt{\det\mathbb G(l_I,\sqrt{-1}\theta_J)}}.
\end{equation}
\end{lemma}

\begin{proof}[Proof of Theorem \ref{main3}] We recall that for each $\epsilon_I\in\{1,-1\}^I,$ $\alpha^*_i=\pi+\epsilon_i\mu_i\sqrt{-1}l_i,$ $\beta_i=\pi+{\mu_i\theta_i}.$ Then by (\ref{CDF}), Lemmas \ref{L5.2}, \ref{L5.3} and \ref{L5.4} and the fact that
$\{1\}=2\sqrt{-1}\sin\frac{2\pi}{r}=\frac{4\pi \sqrt{-1}}{r}+O\Big(\frac{1}{r^2}\Big),$ 
we have
\begin{equation*}
\begin{split}
&\frac{C^{\epsilon_I}}{\sqrt{-\det\mathrm{Hess}\frac{\mathcal W^{\epsilon_I}(z^{\epsilon_I})}{4\pi\sqrt{-1}}}}\\
=&\frac{(-1)^{\frac{|I|+2}{2}+\frac{r}{2}(|I|-2)}n(a_J)}{2^{\frac{7 |I|+2}{2}}\pi^{\frac{3 |I|-2}{2}}} \frac{e^{-\sum_{k=1}^{ 6}\mu_kl_k}}{\big(4\pi\sqrt{-1}\big)^{-\frac{|I|+2}{2}}\sqrt{-(-1)^{\frac{3|I|}{2}}\det\Big(\frac{\partial \theta_{i_1}}{\partial l_{i_2}}\Big)_{i_1,i_2\in I}\det\mathbb G(l_I,\sqrt{-1}\theta_J)}}  r^{\frac{3 |I|-6}{2}}\Big( 1 + O\Big(\frac{1}{r}\Big)\Big).
\end{split}
\end{equation*}
Therefore,
\begin{equation*}
\begin{split}
\sum_{\epsilon_I\in\{1,-1\}^I}&\frac{C^{\epsilon_I}}{\sqrt{-\det\mathrm{Hess}\frac{\mathcal W^{\epsilon_I}(z^{\epsilon_I})}{4\pi\sqrt{-1}}}}\\
=&\frac{(-1)^{\frac{3}{2}+\frac{r}{2}(|I|-2)}n(a_J)}{2^{\frac{3|I|}{2}-1}\pi^{|I|-2}} \frac{e^{-\sum_{k=1}^{6}\mu_kl_k}}{\sqrt{-\det\Big(\frac{\partial \theta_{i_1}}{\partial l_{i_2}}\Big)_{i_1,i_2\in I}\det\mathbb G(l_I,\sqrt{-1}\theta_J)}}  r^{\frac{3 |I|-6}{2}}\Big( 1 + O\Big(\frac{1}{r}\Big)\Big),
\end{split}
\end{equation*}
and by Theorem \ref{dft},
\begin{equation*}
\begin{split}
\mathrm{\widehat {Y}}_r&(\mathbf b^{(r)}_I; \mathbf a^{(r)}_J)\\
&=\frac{(-1)^{\frac{3}{2}+\frac{r}{2}(|I|-2)}n(a_J)}{2^{\frac{3|I|}{2}-1}\pi^{|I|-2}} \frac{e^{-\sum_{k=1}^{6}\mu_kl_k}}{\sqrt{-\det\Big(\frac{\partial \theta_{i_1}}{\partial l_{i_2}}\Big)_{i_1,i_2\in I}\det\mathbb G(l_I,\sqrt{-1}\theta_J)}}  r^{\frac{3 |I|-6}{2}}e^{\frac{r}{\pi}\mathrm{Vol}(\Delta(\theta_I;\theta_J))} \Big( 1 + O\Big(\frac{1}{r}\Big)\Big).
\end{split}
\end{equation*}
\end{proof}


\noindent
Ka Ho Wong\\
Department of Mathematics\\  Yale University\\
New Haven, CT 06520, USA\\
(kaho.wong@yale.edu)
\\

\noindent
Tian Yang\\
Department of Mathematics\\  Texas A\&M University\\
College Station, TX 77843, USA\\
(tianyang@math.tamu.edu)

\end{document}